\documentclass[a4paper,12pt,oneside,reqno]{amsart}
\usepackage[foot]{amsaddr}
\calclayout

\usepackage[utf8]{inputenc}
\usepackage[T1]{fontenc}

\usepackage{amsmath}
\usepackage{amsfonts}
\usepackage{amssymb}
\usepackage{graphicx}
\usepackage{url}
\usepackage{multirow}
\usepackage{caption}
\usepackage{subcaption}
\usepackage{booktabs}

\newcommand{\Nodes}{{\mathcal N}}
\newcommand{\block}{{b}}
\newcommand{\Exp}{{\mathrm Exp}}
\newcommand{\blockID}{\mathrm{blockID}}
\newcommand{\nodeID}{\mathrm{nodeID}}
\newcommand{\blockDelay}{\mathrm{delay}}
\newcommand{\blocks}{\mathrm{Blocks}}
\newcommand{\visibleBlocks}{\mathrm{visBlocks}}
\newcommand{\visibleRef}{\mathrm{visRef}}

\newcommand{\N}{{\mathbb N}}
\renewcommand{\P}{{\mathbb P}}

\newcommand{\E}{{\mathbb E}}
\newcommand{\F}{{\mathcal F}}

\newcommand{\1}{{\mathbf 1}}

\newcommand{\pool}{\mathrm{pool}}
\newcommand{\delay}{{\delta}}
\newcommand{\refBlock}{\mathrm{Ref}}
\renewcommand{\r}{\mathrm{ref}}

\newtheorem{theo}{Theorem}[section]
\newtheorem{ass}[theo]{Assumption}
\newtheorem{lem}[theo]{Lemma}
\newtheorem{defi}[theo]{Definition}
\newtheorem{prop}[theo]{Proposition}
\newtheorem{cor}[theo]{Corollary}
\newtheorem{rem}[theo]{Remark}

\title{Stability of local tip pool sizes}

\author{Sebastian M\"uller}
\address{Sebastian M\"uller (corresponding author), Aix Marseille Universit\'e, CNRS, Centrale Marseille, I2M - UMR 7373, 13453 Marseille, France}
\address{IOTA Foundation, 10405 Berlin, Germany}
\email{sebastian.muller@univ-amu.fr}

\author{Isabel Amigo}
\address{Isabel Amigo, Alexandre Reiffers-Masson, and Santiago Ruano-Rinc\'on, IMT Atlantique, LabSTICC, UMR CNRS 6285, 29238 Brest, France}
\email{isabel.amigo@imt-atlantique.fr}

\author{Alexandre Reiffers-Masson}
\email{alexandre.reiffers-masson@imt-atlantique.fr}

\author{Santiago Ruano-Rinc\'on}
\email{santiago.ruano-rincon@imt-atlantique.fr}

\date{}

\begin{document}

\begin{abstract}
In directed acyclic graph (DAG)-based distributed ledgers, unreferenced blocks (tips) form the backlog of a distributed queueing system. Each new block creates one tip and attempts to remove up to $k$ existing tips by referencing them. With heterogeneous propagation delays, these service decisions are made from delayed local information, so nodes may disagree on the backlog and some reference attempts are wasted.

We study a continuous-time Poisson model with bounded heterogeneous delays and uniform tip selection. We prove that the embedded tip-configuration chain is irreducible, aperiodic, and positive Harris recurrent, and hence admits a unique stationary regime. The observer and local tip-pool sizes have stationary exponential moments, converge to their stationary limits, and satisfy almost-sure ergodic averages. We also derive a Little-type identity relating the stationary mean observer tip count to the mean time until a typical block is first referenced. Simulations are included as qualitative illustrations of the effects of delay variability and issuance heterogeneity.
\end{abstract}


\maketitle

\section{Introduction}
A fundamental challenge in queueing theory arises when servers must make scheduling decisions based on delayed or incomplete state information. In a classical centralised queue, the scheduler observes the full system state and can assign service optimally. When state information is distributed across multiple servers that communicate with nonzero delay, however, service decisions are necessarily based on stale local views, and some service capacity is wasted: two servers may attempt to serve the same job, or a server may attempt to serve a job that has already departed. A natural question is whether such a system remains stable despite this informational handicap.

This paper studies a concrete and practically motivated instance of this problem. In distributed ledger technologies (DLTs), participants maintain a shared append-only database of \emph{blocks}. Each new block must reference a number of previous blocks, and unreferenced blocks are called \emph{tips}. The tip set is precisely the backlog of the system: each new block creates one new job (tip) and simultaneously attempts to serve up to $k$ existing jobs by referencing them. Because blocks propagate through a peer-to-peer network with heterogeneous delays, different nodes maintain different local views of the current backlog, and service decisions are taken from these delayed local perceptions.

Even the simplest distributed ledger architecture, the blockchain, already exhibits this phenomenon. In Bitcoin~\cite{nakamoto2008bitcoin} and similar protocols, multiple miners compete in parallel to extend a linear chain by appending the next block. Each miner works from its own local view, and under nonzero propagation delays several miners may complete a block concurrently, producing \emph{forks}. The protocol resolves forks by the longest-chain rule~\cite{nakamoto2008bitcoin} or a heaviest-subtree variant~\cite{GHOST}, discarding all but one competing block. From a queueing perspective, the discarded blocks represent wasted service: multiple servers attempted to serve the same job (extending the chain tip) using stale local state, and only one attempt was effective. To limit this waste, blockchain throughput is artificially constrained so that each block propagates before the next one is likely to be created -- a design that sacrifices capacity for consistency.

A more recent class of protocols replaces the chain by a directed acyclic graph (DAG) to encode block dependencies. Protocols such as SPECTRE~\cite{spectre}, Byteball~\cite{byteball}, PHANTOM~\cite{phantom}, Prism~\cite{bagaria2019prism}, Aleph~\cite{gkagol2019aleph}, Narwhal~\cite{Narwhal22}, and IOTA~\cite{2020coordicide,OTV} allow multiple nodes to issue blocks concurrently, each referencing several existing tips. Rather than discarding concurrent work, the DAG retains all blocks; concurrency is no longer a failure mode but part of normal operation. This removes the artificial throughput constraint and promises to mitigate problems inherent to the blockchain design such as mining races~\cite{DEVRIES2021105901}, centralisation~\cite{NBERw29396}, miner extractable value~\cite{MEV}, and negative externalities~\cite{Rosenthal22}. However, the queueing problem of wasted service due to delayed information persists in a different form: instead of producing duplicate chain extensions, nodes may now reference the same tip redundantly or reference a tip that has already been served elsewhere.

Within distributed ledger technologies, DAG-based protocols are now a broad research direction; see the recent SoK~\cite{Raikwar2024SoKDAGConsensus}. A useful distinction is between \emph{structured} DAG architectures, where the DAG mainly supports dissemination together with a separate ordering mechanism (e.g., Narwhal/Tusk~\cite{Narwhal22}, Aleph~\cite{gkagol2019aleph}, Prism~\cite{bagaria2019prism}), and \emph{unstructured} block-DAG protocols, where blocks attach directly by references and the ledger remains a persistent DAG (e.g., IOTA 2.0~\cite{OTV,popov2020coordicide}, PHANTOM/GHOSTDAG~\cite{phantom}, and DAGKnight~\cite{SuttonSompolinsky2022DAGKnight}). We focus on the unstructured setting, which corresponds to the richest queueing structure: every node acts as both customer (creating a new tip) and server (referencing existing tips), and no centralised scheduler exists.

The queueing challenge in this setting is that service decisions are taken from delayed local information rather than from a centralised state.  Consequently, service opportunities may be wasted because different nodes can reference the same tip or a tip that has already been referenced elsewhere. The resulting queue-length process is therefore history-dependent and modulated by a random communication environment. The effective service rate is \emph{endogenous}: it depends on the current level of agreement across local views, which itself depends on the backlog and on recent network activity. This makes stability of the tip pool a genuine queueing problem rather than a purely combinatorial property of the DAG.

In this paper, we formalise the evolution of local numbers of unreferenced blocks (tips) and prove their stability under global Poisson issuance and heterogeneous but bounded propagation delays. Except for \cite{kumar2022effect}, such a model is absent from the previous literature, which largely neglects the divergence of local perceptions induced by heterogeneous delays. In \cite{kumar2022effect}, a related distributed model is studied under deterministic delays, deterministic issuance, and discrete time. Our model works in continuous time with random issuance and random delays and leads to stability, stationarity, and queueing-type performance identities. In the next section, we give an informal description of the model.

\subsection{Informal description}
We consider a network of nodes that manage a distributed database. In cryptocurrency applications, this database is called a ledger, but the model could potentially be applied to other use cases of collaborative databases. The data consist of blocks that contain atomic data in the sense that either the entire block is added to the database or all the information in the block is discarded. 
The distributed ledger is assumed to be built using two fundamental mechanisms:

\textbf{Sharing mechanism:} Each node aims to create new blocks and disseminate them to the other nodes. The concrete communication protocol (direct broadcast, gossip, relay network, etc.)\ is not important for our model; what matters is that it induces random, heterogeneous, and bounded propagation delays between nodes. We capture these delays through a standard partial-synchrony assumption (Assumption~\ref{ass:sync}).

\textbf{Reference mechanism:} The blocks in the database are connected to each other by references, forming a directed acyclic graph. The rule is that each newly created block must refer to up to $k\geq 2$ already existing blocks. This reference structure is the defining feature of DAG-based protocols such as IOTA~\cite{popov2015,OTV}, PHANTOM/GHOSTDAG~\cite{phantom}, DAGKnight~\cite{SuttonSompolinsky2022DAGKnight}, Narwhal~\cite{Narwhal22}, and Aleph~\cite{gkagol2019aleph}. The interpretation of a reference depends on the protocol: it may signal verification, validation, or consensus support. In queueing terms, referencing a block corresponds to serving a job. A block that has not yet been referenced by any subsequent block is called a \emph{tip}; the set of tips is the current backlog of the system. The performance of the system depends critically on how nodes choose which tips to reference, and in particular on the fact that previously referenced blocks should no longer be targeted.

Regarding the reference mechanism, we can note that the delay between nodes has a huge impact on the performance of the reference mechanism. Indeed, it is instructive to consider the extreme case where all nodes have the same perception of the database. This can be the case when the creation of a block is instantaneous, i.e.,~there is no delay between selecting the references of the block and sending it to the neighbours, and all neighbours receive the new blocks without delay. Suppose we start the database with one block (the genesis) and assume that no blocks can be created simultaneously. In that case, there will always be only one tip (non-referenced block), as each block is referenced by precisely one other block. However, this situation changes drastically if there is a delay between the selection of references and the time when all nodes have received the new block. In this case, the blocks are created concurrently, and the blocks can be referenced by more than one other block. Thus, a priori, it is no longer clear whether the system is in a stationary regime or the number of tips explodes. 
In this paper, we propose a mathematical procedure to model the different local tip pools and prove the stability of their sizes under standard synchrony assumptions.

\subsection{Contributions}
This paper has four major contributions:
\begin{enumerate}
    \item We formalize the distributed protocol using a stochastic process with global Poisson issuance (Assumption~\ref{ass:poisson}) and a standard partial-synchrony assumption on network delays. To the best of our knowledge, this constitutes the first continuous-time model capturing the local perceptions of a DAG-based distributed ledger and the underlying peer-to-peer communication.
    \item Our main theoretical result, Theorem~\ref{theo:poisson-stability}, proves stability of the observer and local tip pool sizes under Poisson issuance. The proof combines an additive Foster--Lyapunov drift estimate for an augmented process that tracks both the observer tip count and a recent-activity window count with a long-gap minorization / regeneration argument to obtain positive Harris recurrence. A separate exponential-drift argument based on Hajek's theorem then yields exponential tail and moment bounds.
    \item We derive a Little-type identity for the observer tip pool (Section~\ref{sec:little}) in the stationary continuous-time Poisson regime associated with the embedded invariant law, relating the stationary mean number of tips to the mean time until a Palm-typical block is first referenced.
    \item Finally, we present Monte-Carlo simulations (Section~\ref{sec:simulations}) as qualitative illustrations of how the protocol environment influences both the level and the variability of local perceptions.
\end{enumerate}

\medskip
In our model, we introduce a \emph{perfect observer} -- a fictitious entity that sees all blocks instantaneously with zero delay. Its tip pool size $X_n^{(o)}$ is the natural global backlog functional. It is not identical to any local backlog under delay, but Lemma~\ref{lem:X_nX_nc} shows that local and common tip-pool sizes differ from the observer count by at most a term controlled by the recent-activity count $L_n$. The key technical difficulty is that neither component alone admits a simple negative drift: the backlog drift depends on the level of agreement across local views, which is controlled by $L_n$, while $L_n$ fluctuates independently of the backlog. Theorem~\ref{theo:poisson-stability} shows that an augmented Lyapunov function coupling both quantities nevertheless yields a Foster--Lyapunov drift condition.

\subsection{Related work}\label{sec:relworks}

Several works have investigated performance metrics of distributed ledgers, most of them focusing on blockchain-like architectures; see, e.g.,~\cite{Caixiang2020} and references therein. In contrast, formal mathematical models for DAG-based distributed ledgers remain comparatively limited.

\textbf{Early models and centralized views.}
The first mathematical model of a DAG-based distributed ledger was introduced by Popov~\cite{popov2015}. In that model, block arrivals are driven by a Poisson process and the ledger is observed through a \emph{global} view. Communication delays are homogeneous and deterministic, and newly created blocks can be referenced only after a fixed delay. Under these assumptions, Popov heuristically derives an expression for the expected number of tips, relying on a stationarity assumption that is supported by simulations but not proved. Building on this framework, later works studied non-Poisson arrivals~\cite{li2020direct}, discrete-time stochastic models~\cite{bramas2018}, fluid approximations~\cite{ferraro2018distributed,ferraro2018iota}, and simulation-based performance questions~\cite{kusmierz2019properties,park2019performance}. These works do not capture the divergence of local views induced by peer-to-peer propagation delays.

\medskip
\textbf{Heterogeneous delays.}
More recently, several works have attempted to incorporate delay heterogeneity. The simulator proposed in~\cite{zander2019dagsim} explicitly models transmission delays between nodes in a multi-agent DAG-based ledger, but does not provide a mathematical analysis of the resulting dynamics. From a theoretical standpoint, Penzkofer et al.~\cite{penzkofer2021impact} study heterogeneous delays arising from different block processing times while still assuming that a central node maintains a consistent global ledger state. Thus, the heterogeneity in~\cite{penzkofer2021impact} is not due to propagation delays and does not lead to divergent local views. The work closest in spirit to ours is~\cite{kumar2022effect}, which considers a distributed setting where nodes maintain local copies of the ledger and experience delays. However, that model assumes deterministic and time-invariant delays between nodes, as well as a fixed number of issuances per node at each discrete time step. Our work departs from~\cite{kumar2022effect} by allowing \emph{random, heterogeneous, and asymmetric delays}, together with random issuance times in continuous time, and by deriving stationary performance properties for the resulting local-view system.

\medskip
\textbf{Algorithmic improvements and strategic behavior.}
A separate line of work studies improvements of the tip-selection mechanism itself, for instance to enhance security~\cite{attias2019choose,bramas2021efficient}, fairness~\cite{bu2019g}, or robustness to strategic behavior~\cite{popov2019equilibria}. These studies typically build on centralized or homogeneous-delay models and do not address stability of the tip pool under heterogeneous local views.

\medskip
\textbf{Relation to blockchain performance models.}
Several works in the blockchain literature study the impact of network capacity and communication delays on performance and consistency. For example, Prism~\cite{bagaria2019prism} derives throughput and confirmation-time bounds as functions of network capacity, \cite{hegde2022performance} studies ordered-completion times for block dissemination, and \cite{decker2013information} analyzes how network delays induce forks in blockchains. While such systems can also be represented by stochastic DAGs, their attachment rules are fundamentally different: they enforce an eventual chain structure rather than maintaining a persistent DAG through random attachment. Among blockchain-related works, the closest to our perspective is~\cite{gopalan2020stability}, which studies stability properties under deterministic attachment rules. In contrast, we focus on \emph{random attachment under delayed local information} and on the resulting tip-pool dynamics, which are a central performance metric in DAG-based ledgers.

\section{Notations and setting}
The main notation used throughout the paper is collected in Table~\ref{tab:notation}.
\begin{table}[t]
    \centering
    \begin{tabular}{|l|l|}
    \hline
         Variable & Description \\ \hline  \hline
        $\Nodes:=\{1,\ldots,N\}$ & set of nodes\\
        $\lambda_i\in\mathbb{R}_+$ & block issuance rate of node $i$\\ 
        $\lambda:=\sum_{i=1}^N\lambda_i$ & total block issuance rate\\
$\delay_{i}^{(\block)}(j)$& random variable describing \\
& latency from issuer node $i$ to node $j$ for block $\block$\\
$\Delta_i(j)$ & latency distribution from issuer node $i$ to node $j$\\ 
$\Delta\in \mathbb{R}_+$ &  maximal latency between two nodes  \\ 
\hline
$k$  & number of blocks to be referenced by a new block \\ 
\hline
$\pool_n^{-}(i)$ & pre-issuance tip pool of node $i$ at time $t_n$\\
$\pool_n(i)$ & post-issuance tip pool of node $i$ at time $t_n$\\
$\pool_n^{(c)}$ & post-issuance common tip pool at time $t_n$\\
$\pool_n^{(o)}  $ & post-issuance tips of the perfect observer at time $t_n$\\
$X_n^{(i)}:=\lvert \pool_n(i) \rvert$ & size of the post-issuance tip pool of node $i$\\
$X_n^{(c)}:= \lvert\pool_n^{(c)}\rvert$ & size of the post-issuance common tip pool\\
$X_n^{(o)}:= \lvert \pool_n^{(o)}\rvert $ & size of the post-issuance tip pool of the perfect observer\\
\hline
    \end{tabular}
    \caption{Main notation}
    \label{tab:notation}
\end{table}

\subsection{Peer-to-peer network}

We consider a peer-to-peer network with $N$ nodes, $\Nodes:=\{1,\ldots, N\}$, that create and exchange blocks without a central authority. Because blocks propagate through the network with nonzero latency, different nodes hold different local perceptions of the current state.

We model propagation delays as random, asymmetric, and block-dependent. For a block $\block$ issued by node~$i$, the time until node~$j$ receives it is a random variable $\delay_{i}^{(\block)}(j)\ge 0$. These delays are i.i.d.\ across blocks: for every block $\block$ issued by node~$i$, the delay $\delay_{i}^{(\block)}(j)$ is independently distributed as $\Delta_{i}(j)$. By convention, $\delay_{i}^{(\block)}(i)=0$ for all $i\in\Nodes$ and all blocks $\block$, i.e.\ a node sees its own blocks instantaneously.

Three classical synchrony regimes are distinguished in the distributed-systems literature~\cite{DwLySt:88}. In a \emph{fully synchronous} system, the delay distributions $\Delta_{i}(j)$ are bounded by a known constant that is used in the protocol design. In a \emph{fully asynchronous} system, no finite upper bound on delays is assumed. Full synchrony is too restrictive for real-world networks, where delays fluctuate unpredictably, while full asynchrony precludes many useful protocol guarantees. The \emph{partial synchrony} model strikes a practical middle ground: a finite delay bound exists but is unknown and not exploited by the protocol.

\begin{ass}[Partial synchrony]\label{ass:sync}
There exists some $\Delta<\infty$ such that
\[
\P( \Delta_{i}(j)\leq \Delta)=1, \qquad \forall i,j \in \Nodes.
\] The exact value of $\Delta$ is unknown, and its value is not used in the protocol design.
\end{ass}
This assumption is compatible with a wide range of communication mechanisms (direct broadcast, gossip protocols, relay networks): whichever dissemination scheme is used, it suffices that every block reaches every node within a finite time bound.

\subsection{Block issuance}

The blocks are created or issued by the participating nodes at random times $0<t_1<t_2<\cdots$. For notational convenience we also set $t_0:=0$. We model the global issuance pattern by a homogeneous Poisson process and attach to each issuance the marks that determine the issuing node, the block identity, the propagation delays, and the chosen references.

\begin{ass}[Global Poisson issuance]\label{ass:poisson}
The issuance times $(t_n)_{n\ge1}$ form a homogeneous Poisson process on $\mathbb{R}_+$ with rate $\lambda>0$. At each issuance time $t_n$ a mark
\[
\kappa_n=(\blockID_n,\refBlock_n,\nodeID_n,\blockDelay_n)
\]
is attached, and the marks satisfy:
\begin{itemize}
    \item $(\blockID_n)_{n\ge1}$ are i.i.d.\ $\mathrm{Unif}[0,1]$, independent of the Poisson process (hence distinct almost surely);
    \item $(\nodeID_n)_{n\ge1}$ are i.i.d.\ random variables with values in $\Nodes$, independent of everything else;
    \item conditional on $(t_k,\blockID_k,\nodeID_k)_{k\le n}$, the delay vector $\blockDelay_n$ is independent of the past and satisfies Assumption~\ref{ass:sync};
    \item the reference list $\refBlock_n$ is obtained from the issuer's pre-issuance local tip pool, $\pool_n^{-}(\nodeID_n)$, by sampling up to $k$ distinct tips uniformly without replacement.
\end{itemize}
\end{ass}

To encode heterogeneous issuance propensities, we associate to each node $i\in\Nodes$ a nominal rate parameter $\lambda_i>0$ with $\lambda:=\sum_{i=1}^N\lambda_i$; one may think of $\P(\nodeID_n=i)=\lambda_i/\lambda$, although the analysis below only uses that $(\nodeID_n)_{n\ge1}$ is i.i.d.

We briefly describe each mark. The block identifier $\blockID_n\sim\mathrm{Unif}[0,1]$ ensures that blocks are almost surely distinct, mirroring the use of cryptographic hashes in practice. The issuer identity $\nodeID_n\in\Nodes$ records which node created the block. The delay vector $\blockDelay_n=(\delay_{\nodeID_n}^{(\blockID_n)}(j))_{j\in \Nodes}$ gives the propagation time to each node and satisfies Assumption~\ref{ass:sync}.

The reference list $\refBlock_n$ requires more care, as it depends on the history. At epoch $t_n$, we first process at every node all block arrivals with arrival time at most $t_n$, obtaining the pre-issuance local tip pool $\pool_n^{-}(i)$ defined below. If node~$i$ issues the new block, it references
\[
m_n(i):=\min\{k,\lvert \pool_n^{-}(i)\rvert\}
\]
distinct tips sampled uniformly without replacement from that pre-issuance pool. The reference list is written as $\refBlock_n = (\r_1, \ldots, \r_{m_n(\nodeID_n)})$, where each $\r_j$ is the identifier of a previous block, possibly the genesis block $\star$. We denote $(\Omega, \F, \P)$ the underlying probability space and let $\F_n =\sigma ((t_1, \kappa_1),\ldots, (t_n,\kappa_n))$ be the filtration generated by the marked point process. Unlike the other marks, $\refBlock_n$ is not independent of $\F_{n-1}$.

\subsection{Tip selection and dynamics}
At the issuance epochs $t_1,t_2,\ldots$, nodes may have different local views. The set of blocks created up to time $t_n$ is
\begin{equation}
    \blocks_n := \{\star\}\cup\{\blockID_1,\ldots,\blockID_n\},
\end{equation}
where $\star$ denotes the genesis block and we use the convention $t_0:=0$.
The set of blocks created between $t_\ell$ and $t_m$ (with $1\le\ell\le m\le n$) is denoted by 
\begin{equation}
    \blocks_{\ell,m} := \{\blockID_k: \ell\le k\le m\}.
\end{equation}
Due to communication delays, these blocks are not immediately visible to all nodes. At each issuance epoch $t_n$ we distinguish between a pre-issuance state and a post-issuance state. By convention, all messages whose arrival time is at most $t_n$ are processed before the $n$th block is issued. The pre-issuance state is the one used to choose the references of block $n$.

For every node $i$, we therefore define the set of blocks visible just before the issuance at time $t_n$ by
\begin{equation}
    \visibleBlocks_n^{-}(i) := \{\star\}\cup\{\blockID_k: k<n,\ t_k + \blockDelay_k(i) \le t_n\}
\end{equation}
and the corresponding set of visible references by
\begin{equation}
    \visibleRef_n^{-}(i) := \bigcup_{k<n:\ t_k + \blockDelay_k(i) \le t_n} \{\text{entries of }\refBlock_k\},
\end{equation}
where in the last display we regard each reference list $\refBlock_j$ as a finite set of block identifiers. The pre-issuance local tip pool of node $i$ is then
\begin{equation}
    \pool_n^{-}(i) := \visibleBlocks_n^{-}(i) \setminus \visibleRef_n^{-}(i).
\end{equation}
The reference list $\refBlock_n$ is sampled from $\pool_n^{-}(\nodeID_n)$.

Immediately after the issuance of block $n$, we define the post-issuance visible block set of node $i$ by
\begin{equation}
    \visibleBlocks_n(i) := \{\star\}\cup\{\blockID_k: k\le n,\ t_k + \blockDelay_k(i) \le t_n\}
\end{equation}
and the corresponding post-issuance visible references by 
\begin{equation}
    \visibleRef_n(i) := \bigcup_{k\le n:\ t_k + \blockDelay_k(i) \le t_n} \{\text{entries of }\refBlock_k\}.
\end{equation}

\begin{defi}[Different tip pools]
The \emph{post-issuance local tip pool} from node $i\in \Nodes$ at time $t_n$ is defined as
\begin{equation}
    \pool_n(i):=\visibleBlocks_n(i) \setminus \visibleRef_n(i).
\end{equation}
The \emph{post-issuance common tip pool} at time $t_n$ is defined as
\begin{equation}
    \pool_n^{(c)} := \bigcap_{i\in \Nodes} \pool_n(i).
\end{equation}
The \emph{post-issuance (perfect) observer tip pool} at time $t_n$ is defined as
\begin{equation}
    \pool_n^{(o)} := \blocks_n \setminus \bigcup_{k=1}^{n}\{\text{entries of }\refBlock_k\}.
\end{equation}
\end{defi}

\begin{defi}[Tip pool sizes]
We denote by $X_n^{(i)}:= \lvert \pool_n(i) \rvert$ the number of post-issuance tips at node $i$ at time $t_n$. We also define the post-issuance common tip pool size $X_n^{(c)}= \lvert\pool_n^{(c)}\rvert$.
We denote by $X_n^{(o)}= \lvert \pool_n^{(o)}\rvert $ the number of post-issuance tips of the perfect observer.
\end{defi}

The process starts at time $n=0$ with one tip, the genesis block $\star$. More precisely, for every $i\in\Nodes$ we set 
\begin{equation}
    \pool_{0}^{-}(i)=\pool_{0}(i)=\pool_{0}^{(c)}=\pool_{0}^{(o)}=\{\star\},
\end{equation}
and the corresponding tip-pool sizes are equal to $1$.

The pre- and post-issuance tip pool sizes can be defined for all positive real times and can be seen as continuous-time stochastic processes. Due to the delay, the local and common tip pool sizes may even change at times different from the ones given by the point process. However, since nodes issue blocks only at times $t_1,t_2,\ldots$ we observe these processes at those epochs.

Since the observer has zero delays and perceives the blocks right after their creation, it always sees at least as many referenced blocks as any individual node. Hence, when a node issues a block, the tips it references from its pre-issuance local tip pool are immediately deleted from the observer tip pool and the newly issued block is added to the observer tip pool. The newly referenced blocks are also removed, immediately, from the common tip pool, but the new block is added to the common tip pool only after all nodes receive it. If a block had been referenced by some earlier block, then at least the issuer of that earlier block would see the reference immediately, so the block could not belong to every local tip pool. Hence every common tip is an observer tip, and consequently $X_n^{(c)}\leq X_n^{(o)}$ for all $n$.

A crucial observation is that we also have a lower estimate conditioned on the number of blocks recently issued. For $n\ge 1$ let
\[
   L_n := \lvert\{m\in\{1,\ldots,n\}: t_m\in(t_n-\Delta,t_n]\}\rvert,
   \qquad
   L_0:=0.
\]
Then $L_n$ denotes the number of blocks issued in the time interval $(t_n-\Delta,t_n]$. It can also be interpreted as the number of all possibly non-visible blocks at time $t_n$. This definition of $L_n$ also implies that the references chosen at epoch $t_n$ can only depend on the observer tips already present at time $t_{n-L_n}$ together with the new blocks issued between $t_{n-L_n}$ and $t_n$.

\begin{lem}\label{lem:X_nX_nc}
For every $n\ge0$ and every $L\in\{0,\ldots,n\}$, on the event $\{L_n=L\}$ we have, pathwise,
\begin{equation}
     X_n^{(c)}\geq X_n^{(o)}-(k+1)L,
\end{equation}
and
\begin{equation}
     X_n^{(i)}\leq X_n^{(o)}+kL,\qquad \forall i\in\Nodes.
\end{equation}
Equivalently, these inequalities hold almost surely under the conditioning $L_n=L$ whenever $\P(L_n=L)>0$.
\end{lem}
\begin{proof}
	Fix $n\ge0$ and $L\in\{0,\ldots,n\}$ and work on the event $\{L_n=L\}$. By definition of $L_n$, exactly $L$ blocks are issued in the time interval $(t_n-\Delta,t_n]$, namely at times $t_{n-L+1},\ldots,t_n$.
	
	For the transition caused by the issuance at $t_{m+1}$, a single new block is created and becomes a tip in the observer pool, while at most $k$ previously existing observer tips are removed (those referenced by the new block). Hence the observer tip count satisfies, pathwise,
	\[
	   X_{m+1}^{(o)} \le X_m^{(o)} + 1,\qquad X_{m+1}^{(o)} \ge X_m^{(o)} - (k-1),\qquad m\ge 0.
	\]
	Applying these inequalities along the $L$ issuance times in $(t_n-\Delta,t_n]$ yields
	\begin{equation}\label{eq:Xo-step-bounds}
	   X_n^{(o)} \le X_{n-L}^{(o)} + L
	   \quad\text{and}\quad
	   X_n^{(o)} \ge X_{n-L}^{(o)} - (k-1)L.
	\end{equation}
	
	We now relate $X_n^{(c)}$ and $X_n^{(i)}$ to the observer process. First suppose that $n-L\ge1$. Since $L_n=L$, we have $t_{n-L}\le t_n-\Delta$. Thus any block that is a tip in the observer pool at time $t_{n-L}$ has been issued at or before $t_{n-L}$ and therefore has had at least $\Delta$ units of time to propagate by time $t_n$. By Assumption~\ref{ass:sync}, each such block is visible at every node at time $t_n$. If $n-L=0$, the only possible observer tip at time $t_0$ is the genesis block, which is visible at every node by the initial condition. Hence in both cases every tip in $\pool_{n-L}^{(o)}$ is visible at every node at time $t_n$.
	
	Among the $X_{n-L}^{(o)}$ observer tips at time $t_{n-L}$, at most $kL$ distinct tips can be referenced by the $L$ blocks issued in the interval $(t_{n-L},t_n]$, since each of these blocks references at most $k$ tips. Any tip from $\pool_{n-L}^{(o)}$ that is not referenced in $(t_{n-L},t_n]$ therefore remains a tip at every node at time $t_n$: it is older than $\Delta$ (hence visible everywhere) and is not the target of any visible reference. Hence at least $X_{n-L}^{(o)}-kL$ tips are present in \emph{every} local tip pool at time $t_n$, so
	\[
	   X_n^{(c)} \;\ge\; X_{n-L}^{(o)} - kL.
	\]
	Combining this with the first inequality in \eqref{eq:Xo-step-bounds} gives
	\[
	   X_n^{(c)} \;\ge\; X_{n-L}^{(o)} - kL \;\ge\; X_n^{(o)} - (k+1)L,
	\]
	which proves the first claimed inequality.
	
	For the second claim, note that any tip counted by node $i$ but not by the observer can only arise from a reference that the observer has already ``seen'' but node $i$ has not yet seen. By Assumption~\ref{ass:sync}, any block issued at or before time $t_n-\Delta$ is visible at node $i$ by time $t_n$, so only blocks issued in $(t_n-\Delta,t_n]$ can have references that are not yet visible to node $i$. There are exactly $L$ such blocks, and each of them references at most $k$ tips, so there are at most $kL$ referenced tips that are already removed from the observer pool but may still appear as tips at node $i$ at time $t_n$. Therefore node $i$ can have at most $kL$ more tips than the observer:
	\[
	   X_n^{(i)} \;\le\; X_n^{(o)} + kL,\qquad \forall i\in\Nodes.
	\]
	Since all these arguments are pathwise on $\{L_n=L\}$, the lemma follows.
	\end{proof}

\section{Stability under Poisson block issuance}\label{sec:poisson}

In this section we analyze the model under Assumptions~\ref{ass:sync} and~\ref{ass:poisson}, and throughout this section we assume $k\ge 2$. In contrast to the hard-core variant discussed in Section~\ref{sec:hardcore-variant}, the sliding-window count $L_n$ is no longer uniformly bounded, so the stability argument must be adapted. We do this by working with an augmented Lyapunov function involving both the observer tip count $X_n^{(o)}$ and the window count $L_n$.

\subsection{Drift of the window count and of the observer tip count}

We first record a simple drift estimate for the sliding-window count $L_n$.

\begin{lem}[Window-count drift under Poisson issuance]\label{lem:poisson-window-drift}
Let $S_{n+1}:=t_{n+1}-t_n$ denote the inter-arrival time between the $n$th and $(n+1)$st block. Under Assumption~\ref{ass:poisson} the random variables $S_{n+1}$ are i.i.d.\ with $S_{n+1}\sim\Exp(\lambda)$ and independent of $\F_n$. In particular, with
\[
   q := \P(S_{n+1}> \Delta) = e^{-\lambda \Delta}>0,
\]
we have, for every $n\in\N$,
\begin{equation}\label{eq:L-drift-poisson}
   \E\bigl[L_{n+1}-L_n \,\big\vert\, \F_n\bigr] \;\le\; 1 - q\,L_n.
\end{equation}
\end{lem}

\begin{proof}
Fix $n$ and condition on $\F_n$. On the event $\{S_{n+1}>\Delta\}$ there is no issuance in $(t_n,t_n+\Delta]$, hence
\[
   (t_{n+1}-\Delta,t_{n+1}] = (t_{n+1}-\Delta,t_n] \cup \{t_{n+1}\}
\]
contains only the new block at $t_{n+1}$, so $L_{n+1}=1$.

On the complementary event $\{S_{n+1}\le\Delta\}$ the new window $(t_{n+1}-\Delta,t_{n+1}]$ is contained in $(t_n-\Delta,t_{n+1}]$ and differs from $(t_n-\Delta,t_n]$ only by the block at $t_{n+1}$ and by the possible removal of blocks that fall into $(t_n-\Delta,t_{n+1}-\Delta]$. Thus at most one new block enters the window (at $t_{n+1}$), while some of the old-window blocks may leave. In particular,
the number of blocks in the new window cannot exceed the number in the old window plus this single new block, so $L_{n+1}\le L_n+1$. Combining these two inequalities yields
\[
   \E[L_{n+1}\mid\F_n]
   \;\le\; 1\cdot q + (L_n+1)\cdot(1-q)
   \;=\; (1-q)L_n + 1.
\]
\end{proof}

To analyse the next observer increment, we must work with the pre-issuance tip pools at time $t_{n+1}$. We therefore define
\[
   \pool_{n+1}^{-(c)} := \bigcap_{i\in\Nodes}\pool_{n+1}^{-}(i),
   \qquad
   X_{n+1}^{(i,-)} := \lvert\pool_{n+1}^{-}(i)\rvert,
\]
\[
   X_{n+1}^{(c,-)} := \lvert\pool_{n+1}^{-(c)}\rvert.
\]
We also set
\[
   \widehat L_n := \#\{m\le n: t_m\in (t_{n+1}-\Delta,t_n]\}.
\]
Since there are no issuance times in $(t_n,t_{n+1})$, we have $\widehat L_n\le L_n$ almost surely. Moreover, the perfect-observer tip pool just before time $t_{n+1}$ is still $\pool_n^{(o)}$, because no new block is issued between $t_n$ and $t_{n+1}$. Repeating the pathwise argument of Lemma~\ref{lem:X_nX_nc} with the pre-issuance configuration at time $t_{n+1}$ yields
\begin{align}
   X_{n+1}^{(c,-)} &\;\ge\; X_n^{(o)}-(k+1)\widehat L_n \;\ge\; X_n^{(o)}-(k+1)L_n,
   \nonumber\\
   X_{n+1}^{(i,-)} &\;\le\; X_n^{(o)}+k\widehat L_n \;\le\; X_n^{(o)}+kL_n,
   \quad \forall i\in\Nodes. \label{eq:XcXi-poisson}
\end{align}
Let $I_{n+1}:=\nodeID_{n+1}$ be the next issuer and let
$\mathcal{G}_n:=\sigma(\F_n,t_{n+1},I_{n+1})$.
Conditional on $\mathcal{G}_n$, the $(n+1)$st block samples its references uniformly without replacement from $\pool_{n+1}^{-}(I_{n+1})$. Denote by $R_{n+1}$ the event that the first two sampled references belong to the common pool $\pool_{n+1}^{-(c)}$.

\begin{lem}[Common-reference probability]\label{lem:common-ref-prob}
Define
\[
   \gamma(\phi) := \frac{1-(k+1)\phi}{1+k\phi},\qquad \phi\in[0,(k+1)^{-1}).
\]
On the event $\{X_n^{(o)}=x,\,L_n=\ell\}$ with $x-(k+1)\ell\ge 2$ and $x>0$, writing $\phi:=\ell/x$,
\begin{equation}\label{eq:p-x-ell}
   \P(R_{n+1}\mid \mathcal{G}_n)
   \;\ge\; \gamma(\phi)^2 - \frac{\gamma(\phi)}{x(1+k\phi)}.
\end{equation}
The function $\gamma$ is continuous and strictly decreasing with $\gamma(0)=1$. Hence we can choose $\beta\in(0,(k+1)^{-1})$ sufficiently small that
\begin{equation}\label{eq:gamma-beta}
   \gamma(\beta)^2 > \tfrac34.
\end{equation}
\end{lem}

\begin{proof}
We realise the uniform sample without replacement sequentially; this is only a convenient representation of the unordered reference list. On the event $\{X_{n+1}^{(I_{n+1},-)}\ge2\}$, the probability that the first two sampled references are both common tips equals
\[
   \P(R_{n+1}\mid \mathcal{G}_n)
   = \frac{X_{n+1}^{(c,-)}}{X_{n+1}^{(I_{n+1},-)}}
     \cdot \frac{X_{n+1}^{(c,-)}-1}{X_{n+1}^{(I_{n+1},-)}-1}.
\]
Since $X_{n+1}^{(I_{n+1},-)}-1\le X_{n+1}^{(I_{n+1},-)}$, we have
\[
   \frac{X_{n+1}^{(c,-)}-1}{X_{n+1}^{(I_{n+1},-)}-1}
   \ge \frac{X_{n+1}^{(c,-)}-1}{X_{n+1}^{(I_{n+1},-)}},
\]
so it suffices to bound the two ratios with denominator $X_{n+1}^{(I_{n+1},-)}$. From \eqref{eq:XcXi-poisson} and the definition of $\phi$ we have
\[
   \frac{X_{n+1}^{(c,-)}}{X_{n+1}^{(I_{n+1},-)}}
   \ge \frac{x-(k+1)\ell}{x+k\ell}
   = \frac{1-(k+1)\phi}{1+k\phi}
   = \gamma(\phi).
\]
Similarly,
\[
   \frac{X_{n+1}^{(c,-)}-1}{X_{n+1}^{(I_{n+1},-)}}
   \ge \frac{x-(k+1)\ell-1}{x+k\ell}
   = \gamma(\phi) - \frac{1}{x+k\ell}.
\]
Multiplying these two bounds and using $x+k\ell = x(1+k\phi)$ gives \eqref{eq:p-x-ell}.
\end{proof}

\begin{lem}[Conditional drift of the observer tip count]\label{lem:X-drift-poisson}
There exist constants $\beta\in(0,(k+1)^{-1})$, $K_0<\infty$ and $\varepsilon_X>0$ such that for all $n\in\N$,
\begin{equation}\label{eq:X-drift-poisson}
   \E\bigl[X_{n+1}^{(o)}-X_n^{(o)} \,\big\vert\, \F_n\bigr]
   \;\le\; -\varepsilon_X
\end{equation}
on the event $\{X_n^{(o)}\ge K_0,\,L_n\le \beta X_n^{(o)}\}$.
\end{lem}

\begin{proof}
Fix $\beta$ as above and choose $K_0$ so large that
\[
   \frac{1}{x} \le \gamma(\beta)^2 - \tfrac34,
   \qquad
   x\bigl(1-(k+1)\beta\bigr)\ge 2,
   \qquad \forall x\ge K_0.
\]
Fix $n\in\N$ and work on the event $\{X_n^{(o)}=x,\,L_n=\ell\}$ with $x\ge K_0$ and $\ell\le\beta x$. Then $0\le\phi=\ell/x\le\beta<(k+1)^{-1}$. Since $\gamma$ is decreasing in $\phi$, this implies $\gamma(\phi)\ge\gamma(\beta)$ and hence $\gamma(\phi)^2\ge\gamma(\beta)^2$. Moreover, for such $\phi$ we have
\[
   0\le \frac{\gamma(\phi)}{1+k\phi} = \frac{1-(k+1)\phi}{(1+k\phi)^2} \le 1,
\]
so \eqref{eq:p-x-ell} yields
\[
   \P(R_{n+1}\mid \mathcal{G}_n)
   \;\ge\; \gamma(\phi)^2 - \frac{\gamma(\phi)}{x(1+k\phi)}
   \;\ge\; \gamma(\beta)^2 - \frac{1}{x}
   \;\ge\; \tfrac34.
\]
The second condition in the choice of $K_0$ ensures that $X_{n+1}^{(c,-)}\ge 2$ on this region, so the event $R_{n+1}$ is well defined. A common pre-issuance tip cannot have been referenced by any earlier block, because the issuer of such a referencing block would already have removed it from its own local tip pool. Thus every tip in $\pool_{n+1}^{-(c)}$ is an observer tip just before time $t_{n+1}$. Whenever $R_{n+1}$ occurs, at least two observer tips are removed and at most one new tip is created (since sampling is without replacement), so $X_{n+1}^{(o)}-X_n^{(o)}\le -1$ on that event. On the complementary event the increment is bounded above by $+1$, since each issuance creates at most one new tip. Therefore,
\[
   \E[X_{n+1}^{(o)}-X_n^{(o)}\mid \mathcal{G}_n]
   \;\le\; (-1)\cdot \P(R_{n+1}\mid \mathcal{G}_n) + 1\cdot\bigl(1-\P(R_{n+1}\mid \mathcal{G}_n)\bigr)
   \;\le\; 1-2\cdot\tfrac34 = -\tfrac12
\]
on the event $\{X_n^{(o)}\ge K_0,\,L_n\le\beta X_n^{(o)}\}$. Taking conditional expectation once more with respect to $\F_n$ proves \eqref{eq:X-drift-poisson} with $\varepsilon_X:=1/2$.
\end{proof}

\subsection{Lyapunov drift for the augmented process}

We now define the Lyapunov function
\[
   V_n := X_n^{(o)} + \alpha L_n,
\]
where $\alpha>0$ is a constant to be chosen. Note that $V_n$ is nonnegative and $\F_n$-measurable for each $n$.

\begin{prop}[Long-gap minorization and small sets]\label{prop:small-set-poisson}
For a fixed choice of $\alpha>0$, let $\mathsf{Y}_n$ denote the augmented embedded state consisting of the observer tip pool, all local tip pools, and the marked configuration of blocks issued in the last delay window:
\[
   \mathsf{Y}_n
   :=
   \Bigl(
      \pool_n^{(o)},
      (\pool_n(i))_{i\in\Nodes},
      \bigl((t_n-t_m,\kappa_m)\bigr)_{m\le n:\, t_m\in(t_n-\Delta,t_n]}
   \Bigr).
\]
Then $(\mathsf{Y}_n)_{n\ge0}$ is a time-homogeneous Markov chain on a standard Borel state space. Moreover, there exists a probability measure $\nu_\ast$ such that for every $M<\infty$, with
\[
   C_M:=\{y:\ V(y)\le M\},
\]
there are integers $m_M\in\N$ and constants $p_M,\widetilde p_M>0$ for which
\[
   P^{m_M+2}(y,\cdot)\ge p_M\,\nu_\ast(\cdot),
   \qquad
   P^{m_M+3}(y,\cdot)\ge \widetilde p_M\,\nu_\ast(\cdot),
   \qquad \forall y\in C_M,
\]
where $P$ is the transition kernel of $(\mathsf{Y}_n)$. In particular, each $C_M$ is a small set and the chain is $\nu_\ast$-irreducible and aperiodic.
\end{prop}

\begin{proof}
Only blocks issued during the most recent interval of length $\Delta$ can still arrive at different nodes in the future. The current local pools record all tips that may still be selected by some node, while the recent marked configuration determines which not-yet-globally-visible blocks and references arrive before the next issuance. Thus the enlarged state determines all pre-issuance local tip pools at the next decision epoch and therefore determines the law of the next state. By the memoryless property of exponential inter-arrival times and the independence of future marks, $(\mathsf{Y}_n)$ is a time-homogeneous Markov chain.

To make the state space explicit, let $\widehat{\mathsf B}:=\{\star\}\sqcup [0,1]$, where $\star$ denotes the genesis block, and let
\[
   \mathsf K
   :=
   \widehat{\mathsf B}\times
   \Bigl(\bigsqcup_{r=0}^k \widehat{\mathsf B}^{\,r}\Bigr)\times
   \Nodes\times [0,\Delta]^N
\]
be the mark space. Equip $\widehat{\mathsf B}$ with any Borel total order extending $\star<x$ for all $x\in[0,1]$, and for $\ell\ge 1$ write
\[
   \widehat{\mathsf B}^{\,\ell}_{\uparrow}
   :=
   \{(b_1,\ldots,b_\ell)\in \widehat{\mathsf B}^{\,\ell}: b_1<\cdots < b_\ell\}.
\]
We use the convention that $\widehat{\mathsf B}^{\,0}_{\uparrow}$ is a singleton.
Encoding each finite tip pool by the increasing tuple of its elements, the state space of $(\mathsf{Y}_n)$ can be identified with the countable disjoint union
\[
   \mathsf Y
   :=
   \bigsqcup_{\ell_o\ge 1}
   \bigsqcup_{\ell_1,\ldots,\ell_N\ge 0}
   \bigsqcup_{m\ge 0}
   \widehat{\mathsf B}^{\,\ell_o}_{\uparrow}
   \times
   \prod_{i=1}^N \widehat{\mathsf B}^{\,\ell_i}_{\uparrow}
   \times \bigl([0,\Delta]\times \mathsf K\bigr)^m,
\]
which is a standard Borel space because each component is a Borel subset of a Polish space and a countable disjoint union of standard Borel spaces is standard Borel.
We take the actual state space to be the measurable subset of $\mathsf Y$ consisting of admissible configurations generated by the protocol dynamics; this subset is again standard Borel.

Let $\nu_\ast$ be the law of the enlarged post-issuance state obtained when one starts from a synchronised configuration with a single common tip and no recent blocks, and then performs two fresh issuances with inter-arrival gaps larger than $\Delta$. This law includes the observer pool, all local pools, and the remaining recent marked configuration after the second issuance. Fix $M<\infty$ and set $m_M:=\max\{1,\lceil M\rceil\}$. If $\mathsf{Y}_n=y\in C_M$, then necessarily $X_n^{(o)}\le M$. Consider the events
\begin{align*}
   E_M&:=\{t_{n+j}-t_{n+j-1}>\Delta,\ \forall j=1,\dots,m_M+2\},\\
   \widetilde E_M&:=\{t_{n+j}-t_{n+j-1}>\Delta,\ \forall j=1,\dots,m_M+3\}.
\end{align*}
Since the inter-arrival gaps are i.i.d.\ exponential with rate $\lambda$, these events have probabilities
\[
   p_M = \P(E_M)=e^{-\lambda\Delta (m_M+2)}>0,
   \qquad
   \widetilde p_M = \P(\widetilde E_M)=e^{-\lambda\Delta (m_M+3)}>0.
\]
On either event, before each of the next issuances all previously issued blocks are visible at every node, so the recent window has emptied and the local, common, and observer tip pools coincide at each decision epoch. Because $k\ge 2$, whenever the current synchronised observer tip count is larger than $1$, the next synchronised issuance removes at least two distinct tips and creates one new tip, so the observer tip count decreases by at least $1$. Once the observer tip count reaches $1$, further synchronised issuances keep it equal to $1$. As $X_n^{(o)}\le M\le m_M$, after at most $m_M$ such synchronised issuances the chain is in a post-issuance one-tip configuration whose unique tip was created during the long-gap event itself and is therefore independent of the initial state $y$. The next two synchronised issuances, still on $E_M$ or $\widetilde E_M$, then produce an enlarged post-state whose law depends only on the fresh marks generated during those last two issuances; by definition, that law is $\nu_\ast$. Therefore
\[
   P^{m_M+2}(y,\cdot)\ge p_M\,\nu_\ast(\cdot),
   \qquad
   P^{m_M+3}(y,\cdot)\ge \widetilde p_M\,\nu_\ast(\cdot),
   \qquad \forall y\in C_M,
\]
so $C_M$ is small.

For $\nu_\ast$-irreducibility, fix an arbitrary initial state $y$ and write $x:=X^{(o)}(y)<\infty$. On the event of $x+2$ consecutive inter-arrival gaps larger than $\Delta$, the same synchronisation argument shows that after at most $x$ such issuances the observer tip count has reached $1$, and that the final two long-gap issuances produce a state with law $\nu_\ast$. Hence for every measurable set $A$,
\[
   P^{x+2}(y,A)\ge e^{-\lambda\Delta(x+2)}\,\nu_\ast(A).
\]
Therefore the chain is $\nu_\ast$-irreducible.

Finally, for any $M\ge 1+\alpha$, the construction of $\nu_\ast$ gives $\nu_\ast(C_M)=1$. The two minorization relations on such a $C_M$ occur at the coprime times $m_M+2$ and $m_M+3$, so $\gcd(m_M+2,m_M+3)=1$. This is the standard strong aperiodicity criterion for an irreducible chain, hence $(\mathsf{Y}_n)$ is aperiodic.
\end{proof}

\begin{theo}[Stability under Poisson issuance]\label{theo:poisson-stability}
Suppose that Assumptions~\ref{ass:sync} and~\ref{ass:poisson} hold and that $k\ge2$. Then there exist constants $\alpha>0$, $K<\infty$ and $\varepsilon>0$ such that
\begin{equation}\label{eq:V-drift-poisson}
   \E\bigl[V_{n+1}-V_n \,\big\vert\, \F_n\bigr]
   \;\le\; -\varepsilon
   \qquad\text{on the event } \{V_n\ge K\},\ \forall n\ge 0.
\end{equation}
Moreover, for the Markov chain $(\mathsf{Y}_n)$ of Proposition~\ref{prop:small-set-poisson}, there exists $b<\infty$ such that
\[
   PV(y)\le V(y)-\varepsilon + b\,\mathbf 1_{C_K}(y),
   \qquad
   C_K:=\{y:\ V(y)\le K\}.
\]
Consequently, $(\mathsf{Y}_n)$ is $\nu_\ast$-irreducible, aperiodic, and positive Harris recurrent. In particular, it admits a unique invariant probability measure $\pi$.
\end{theo}

\begin{proof}
Fix $\beta$, $K_0$ and $\varepsilon_X$ as in Lemma~\ref{lem:X-drift-poisson}. On the region
\[
   A:=\{X_n^{(o)}\ge K_0,\ L_n\le\beta X_n^{(o)}\}
\]
we have
\begin{equation}\label{eq:region-A}
   \E[X_{n+1}^{(o)}-X_n^{(o)}\mid\F_n] \le -\varepsilon_X.
\end{equation}

For all states we have the trivial bound
\begin{equation}\label{eq:X-increment-bound}
   X_{n+1}^{(o)}-X_n^{(o)} \le 1,
\end{equation}
since each new block creates at most one tip. Combining \eqref{eq:L-drift-poisson} and \eqref{eq:X-increment-bound} we obtain, for all $n$,
\begin{equation}\label{eq:V-drift-general}
   \E[V_{n+1}-V_n \mid \F_n]
   \;\le\; \E[X_{n+1}^{(o)}-X_n^{(o)}\mid\F_n] + \alpha\,(1-qL_n)
   \;\le\; 1 + \alpha - \alpha q L_n.
\end{equation}

We now distinguish two regions in the $(x,\ell)$--plane.

\smallskip
\noindent\emph{Region A:} $\{x\ge K_0,\ \ell\le\beta x\}$. On this region we combine \eqref{eq:region-A} and the inequality $\E[L_{n+1}-L_n\mid\F_n]\le 1$ to get
\[
   \E[V_{n+1}-V_n\mid\F_n]
   = \E[X_{n+1}^{(o)}-X_n^{(o)}\mid\F_n]
     + \alpha\,\E[L_{n+1}-L_n\mid\F_n]
   \le -\varepsilon_X + \alpha.
\]
Choosing $\alpha\le \varepsilon_X/2$ yields
\begin{equation}\label{eq:V-drift-A}
   \E[V_{n+1}-V_n\mid\F_n]
   \;\le\; -\tfrac12\varepsilon_X
   \qquad\text{on }\{X_n^{(o)}\ge K_0,\ L_n\le\beta X_n^{(o)}\}.
\end{equation}

\smallskip
\noindent\emph{Region B:} $\{\ell\ge \ell_0\}$ for a suitable $\ell_0$. From \eqref{eq:V-drift-general} we have
\[
   \E[V_{n+1}-V_n\mid\F_n]
   \;\le\; 1+\alpha - \alpha q L_n.
\]
Thus, on the event $\{L_n\ge \ell_0\}$ we get
\begin{equation}\label{eq:V-drift-B}
   \E[V_{n+1}-V_n\mid\F_n]
   \;\le\; 1+\alpha - \alpha q \ell_0.
\end{equation}
Choosing
\[
   \ell_0 := \Bigl\lceil \frac{2(1+\alpha)}{\alpha q}\Bigr\rceil
\]
we obtain
\[
   1+\alpha - \alpha q \ell_0 \le 1+\alpha - 2(1+\alpha) = - (1+\alpha) \le -1.
\]

\smallskip
We now choose
\[
   K := \max\{K_0,\ \ell_0/\beta\} + \alpha \ell_0 + 1.
\]
If $L_n<\ell_0$ and $X_n^{(o)}<K_0$, then
\[
   V_n = X_n^{(o)}+\alpha L_n < K_0+\alpha\ell_0 < K.
\]
If $L_n<\ell_0$, $X_n^{(o)}\ge K_0$ and $L_n>\beta X_n^{(o)}$, then $X_n^{(o)}<L_n/\beta<\ell_0/\beta$, and hence
\[
   V_n = X_n^{(o)}+\alpha L_n < \ell_0/\beta + \alpha \ell_0 < K.
\]
Therefore, whenever $V_n\ge K$ we must have either $L_n\ge \ell_0$ or $X_n^{(o)}\ge K_0$ together with $L_n\le\beta X_n^{(o)}$. Hence the event $\{V_n\ge K\}$ is contained in the union of the two regions
\[
   A := \{X_n^{(o)}\ge K_0,\ L_n\le\beta X_n^{(o)}\},\qquad
   B := \{L_n\ge \ell_0\}.
\]
On $B$ we use the drift estimate \eqref{eq:V-drift-B}, while on $A\setminus B$ we use \eqref{eq:V-drift-A}. Therefore,
\[
   \E[V_{n+1}-V_n\mid\F_n]
   \;\le\; -\varepsilon,
\]
whenever $V_n\ge K$, where we may take
\[
   \varepsilon := \min\left\{\tfrac12\varepsilon_X,\ 1\right\}>0.
\]
This proves \eqref{eq:V-drift-poisson}.

To convert this into a Foster--Lyapunov inequality for the Markov chain $(\mathsf{Y}_n)$, note that
\[
   V_{n+1}-V_n
   =
   (X_{n+1}^{(o)}-X_n^{(o)}) + \alpha(L_{n+1}-L_n)
   \le 1+\alpha
\]
pathwise, because $X_{n+1}^{(o)}-X_n^{(o)}\le 1$ and $L_{n+1}-L_n\le 1$. Hence, writing $PV(y):=\E[V_{n+1}\mid \mathsf{Y}_n=y]$, we obtain
\[
   PV(y)\le V(y)-\varepsilon + b\,\mathbf 1_{C_K}(y),
   \qquad
   b:=1+\alpha+\varepsilon.
\]
Equivalently, with $W:=V/\varepsilon$,
\[
   PW(y)\le W(y)-1+\frac{b}{\varepsilon}\,\mathbf 1_{C_K}(y).
\]

By Proposition~\ref{prop:small-set-poisson}, $(\mathsf{Y}_n)$ is $\nu_\ast$-irreducible and aperiodic, and $C_K$ is a small set, hence also petite. The function $W$ is everywhere finite and bounded on $C_K$. Therefore Meyn and Tweedie~\cite[Theorem~11.3.4]{MeynTweedie2009} applies and yields that $(\mathsf{Y}_n)$ is positive Harris recurrent. A positive Harris recurrent irreducible chain admits a unique invariant probability measure; see Meyn and Tweedie~\cite[Theorems~10.4.4 and~10.4.10(ii)]{MeynTweedie2009}. This completes the proof.
\end{proof}

\subsection{Hitting-time and tail bounds}\label{sec:quantitative}

The additive drift of Theorem~\ref{theo:poisson-stability} yields positive Harris recurrence, but it does not by itself provide the stationary moment bounds used later. To obtain quantitative control we verify Hajek's exponential-drift hypotheses~\cite[Theorem~2.3]{Hajek:82} directly for the same Lyapunov function $V_n$.

\begin{theo}[Hitting-time and tail bounds]\label{theo:hitting}
Under Assumptions~\ref{ass:sync} and~\ref{ass:poisson}, with $k\ge 2$, with $\alpha$ and $K$ chosen as in Theorem~\ref{theo:poisson-stability}, and with $\ell_0$ as defined in its proof, there exist constants $\eta>0$ and $\rho\in(0,1)$ such that for all integers $m,n\ge 0$,
\begin{align}
    \E\bigl[e^{\eta V_{m+n}} \,\big\vert\, \F_m\bigr]
      &\le \rho^n\, e^{\eta V_m}
           + \frac{1-\rho^n}{1-\rho}\,e^{\eta(1+\alpha)}\,e^{\eta K},
      \label{eq:exp-moment}\\[0.4em]
    \P\bigl(V_{m+n}\ge M \,\big\vert\, \F_m\bigr)
      &\le \rho^n\, e^{\eta(V_m - M)}
           + \frac{1-\rho^n}{1-\rho}\,e^{\eta(1+\alpha)}\,e^{\eta(K-M)},
      \label{eq:tail-bound}\\[0.4em]
    \P\bigl(\tau_{K,m}>n \,\big\vert\, \F_m\bigr)
      &\le e^{\eta(V_m - K)}\,\rho^n,
      \label{eq:return-time}
\end{align}
where $\tau_{K,m}:=\min\{j\ge 0: V_{m+j}\le K\}$. In particular, since $X_n^{(o)}\le V_n$ and $X_n^{(i)}\le (1+k/\alpha)\,V_n$ for all $i\in\Nodes$, the observer and local tip pool sizes inherit exponential tail bounds, up to the deterministic rescaling of thresholds.
\end{theo}

\begin{proof}
Set
\[
   \Delta_n:=V_{n+1}-V_n.
\]

We first verify Hajek's upper-region condition on the two regions used in the proof of Theorem~\ref{theo:poisson-stability}.

\smallskip
\noindent\emph{Region A.} On
\[
   A:=\{X_n^{(o)}\ge K_0,\ L_n\le \beta X_n^{(o)}\},
\]
Lemma~\ref{lem:X-drift-poisson} was proved by introducing the event $R_{n+1}$ that the first two references sampled by the next issuer belong to the common pre-issuance tip pool. On $A$ we have
\[
   \P(R_{n+1}\mid \mathcal{G}_n)\ge \tfrac34,
\]
where $\mathcal{G}_n=\sigma(\F_n,t_{n+1},I_{n+1})$. Because $k\ge 2$ and the references are sampled without replacement, the event $R_{n+1}$ implies that two distinct observer tips are removed while one new tip is created, so
\[
   X_{n+1}^{(o)}-X_n^{(o)}\le -1
   \qquad\text{on }R_{n+1}.
\]
Moreover, always $L_{n+1}-L_n\le 1$. Hence
\[
   \Delta_n\le \alpha-1
   \qquad\text{on }R_{n+1},
\]
while on $R_{n+1}^{c}$ we have only the trivial bound
\[
   \Delta_n\le 1+\alpha.
\]
Therefore, on the region $A$,
\[
   \E[e^{\eta\Delta_n}\mid \mathcal{G}_n]
   \le \tfrac34 e^{\eta(\alpha-1)}+\tfrac14 e^{\eta(1+\alpha)}
   =:\rho_A(\eta).
\]
Our choice of $\alpha$ in Theorem~\ref{theo:poisson-stability} satisfies $\alpha\le \varepsilon_X/2=1/4<1/2$, so
\[
   \rho_A(0)=1,
   \qquad
   \rho_A'(0)=\tfrac34(\alpha-1)+\tfrac14(1+\alpha)=\alpha-\tfrac12<0.
\]
By continuity, there exists $\eta_A>0$ such that $\rho_A(\eta)<1$ for all $0<\eta\le \eta_A$. Taking conditional expectation once more with respect to $\F_n$ preserves the same bound on $A$.

\smallskip
\noindent\emph{Region B.} On
\[
   B:=\{L_n\ge \ell_0\},
\]
let $H_n:=\{t_{n+1}-t_n>\Delta\}$. By Assumption~\ref{ass:poisson},
\[
   \P(H_n\mid \F_n)=q=e^{-\lambda\Delta}.
\]
On $H_n$ we have $L_{n+1}=1$, so
\[
   \Delta_n
   =
   (X_{n+1}^{(o)}-X_n^{(o)})+\alpha(1-L_n)
   \le 1+\alpha(1-L_n)
   \le 1+\alpha(1-\ell_0)
\]
throughout the region $B$. On $H_n^c$ we again use only the pathwise bound $\Delta_n\le 1+\alpha$. Hence on $B$,
\[
   \E[e^{\eta\Delta_n}\mid \F_n]
   \le q\,e^{\eta(1+\alpha(1-\ell_0))} + (1-q)\,e^{\eta(1+\alpha)}
   =:\rho_B(\eta).
\]
By the choice of $\ell_0$ in Theorem~\ref{theo:poisson-stability},
\[
   \rho_B(0)=1,
   \qquad
   \rho_B'(0)=q\bigl(1+\alpha(1-\ell_0)\bigr)+(1-q)(1+\alpha)
   = 1+\alpha-\alpha q\ell_0<0.
\]
Therefore there exists $\eta_B>0$ such that $\rho_B(\eta)<1$ for all $0<\eta\le \eta_B$.

Now choose
\[
   0<\eta\le \min\{\eta_A,\eta_B\}
   \qquad\text{and}\qquad
   \rho:=\max\{\rho_A(\eta),\rho_B(\eta)\}<1.
\]
By the proof of Theorem~\ref{theo:poisson-stability}, the event $\{V_n\ge K\}$ is contained in $A\cup B$, and hence
\begin{equation}\label{eq:hajek-A3}
   \E[e^{\eta(V_{n+1}-V_n)}\mid \F_n]\le \rho
   \qquad\text{on }\{V_n\ge K\}.
\end{equation}

We next verify the lower-region condition. If $V_n\le K$, then
\[
   V_{n+1}-K
   = (V_{n+1}-V_n)+(V_n-K)
   \le 1+\alpha
\]
pathwise, because $V_{n+1}-V_n\le 1+\alpha$ and $V_n-K\le 0$. Therefore
\begin{equation}\label{eq:hajek-A4}
   \E[e^{\eta(V_{n+1}-K)}\mid \F_n]\le e^{\eta(1+\alpha)}
   \qquad\text{on }\{V_n\le K\}.
\end{equation}

Fix an integer $m\ge 0$ and apply Hajek's theorem~\cite[Theorem~2.3]{Hajek:82} to the shifted process
\[
   Y_j:=V_{m+j},
   \qquad
   \mathcal{H}_j:=\F_{m+j},
   \qquad j\ge 0.
\]
The bounds \eqref{eq:hajek-A3}--\eqref{eq:hajek-A4} are uniform in the starting time, so Hajek's conclusion gives \eqref{eq:exp-moment}--\eqref{eq:return-time}. The tail bound \eqref{eq:tail-bound} is obtained from \eqref{eq:exp-moment} by Markov's inequality, exactly as in Hajek's theorem.
\end{proof}

Bounds \eqref{eq:tail-bound} and \eqref{eq:return-time} are uniform in the starting time $m$: they depend only on the current value $V_m$, not on the detailed past. From an operational perspective, they show that atypically large tip pools are exponentially unlikely and that returns to a fixed sublevel set have exponentially decaying tails.

A further consequence is a uniform bound on all polynomial moments.

\begin{cor}[Moment bounds]\label{cor:moments}
Under the assumptions of Theorem~\ref{theo:poisson-stability}, for every $r>0$ there exists a constant $c_r<\infty$ such that
\[
   \sup_{n\ge0}\E\bigl[(X_n^{(o)})^r\bigr] \le c_r.
\]
The same holds for $X_n^{(i)}$, $i\in\Nodes$, with constants that can be chosen uniformly over $i$.
\end{cor}

\begin{proof}
The exponential moment bound \eqref{eq:exp-moment} with $m=0$ gives $\sup_{n\ge0}\E[e^{\eta V_n}]<\infty$. For every real $r>0$ there is a finite constant $C_{r,\eta}$ such that $x^r\le C_{r,\eta}e^{\eta x}$ for all $x\ge0$; for instance, one may take $C_{r,\eta}:=\sup_{x\ge0}x^r e^{-\eta x}<\infty$. Hence $\sup_{n\ge0}\E[V_n^r]<\infty$. The claim follows from $X_n^{(o)}\le V_n$ and $X_n^{(i)}\le(1+k/\alpha)V_n$.
\end{proof}

\subsection{Regeneration structure}

The long-gap minorization from Proposition~\ref{prop:small-set-poisson} has a concrete regenerative interpretation. We record it explicitly because it identifies i.i.d.\ post-regeneration cycle segments for the embedded chain.

We define the regeneration indices by
\begin{equation}\label{eq:rho_r}
   \rho_0
   :=
   \inf\{n\ge 2:\ X_{n-2}^{(o)}=1,\ t_{n-1}-t_{n-2}>\Delta,\ t_n-t_{n-1}>\Delta\},
\end{equation}
and, for $r\ge0$,
\begin{equation}\label{eq:rho_r_plus}
   \rho_{r+1}
   :=
   \inf\{n>\rho_r:\ X_{n-2}^{(o)}=1,\ t_{n-1}-t_{n-2}>\Delta,\ t_n-t_{n-1}>\Delta\}.
\end{equation}
At time $t_{\rho_r}$ the chain has just completed two synchronised long-gap issuances starting from a one-tip state. Since $t_{\rho_r}-t_{\rho_r-1}>\Delta$, the only block in $(t_{\rho_r}-\Delta,t_{\rho_r}]$ is the one issued at $t_{\rho_r}$, so $L_{\rho_r}=1$. Moreover, the observer has exactly one tip after the issuance. Hence
\begin{equation}\label{eq:regen-state}
   X_{\rho_r}^{(o)}=1,\qquad L_{\rho_r}=1,\qquad V_{\rho_r}=1+\alpha.
\end{equation}

\begin{prop}[Long-gap regeneration cycles]\label{prop:regeneration-poisson}
Under Assumptions~\ref{ass:sync} and~\ref{ass:poisson}, the indices $(\rho_r)_{r\ge0}$ are almost surely finite stopping times. Moreover, the post-regeneration cycle segments
\[
   \bigl(\mathsf{Y}_{\rho_r+n}\bigr)_{0\le n<\rho_{r+1}-\rho_r},
   \qquad r\ge0,
\]
are i.i.d.
\end{prop}

\begin{proof}
The times $\rho_r$ are stopping times because the events in \eqref{eq:rho_r} and \eqref{eq:rho_r_plus} depend only on the marks up to the corresponding index. On the event $\{n=\rho_r\}$ we have $X_{n-2}^{(o)}=1$ and both gaps $t_{n-1}-t_{n-2}$ and $t_n-t_{n-1}$ exceed $\Delta$. Since there are no issuances in $(t_{n-2},t_{n-1})$ and the unique observer tip present at time $t_{n-2}$ is older than $\Delta$ by time $t_{n-1}-$, Assumption~\ref{ass:sync} implies that all nodes see the same unique tip just before the issuance at time $t_{n-1}$. The first long-gap issuance therefore starts from a synchronised one-tip configuration. The second long gap implies that before the issuance at time $t_n$ all nodes again see the same unique tip and no recent unseen blocks remain. Consequently, the post-$t_n$ state $\mathsf{Y}_n$ depends only on the fresh marks of the issuances at times $t_{n-1}$ and $t_n$, together with the condition that the two gaps exceed $\Delta$, and has the fixed law $\nu_\ast$ from Proposition~\ref{prop:small-set-poisson}. The stopping rule does not otherwise bias these final two fresh marks. Since $(\mathsf{Y}_n)$ is Markov with respect to $(\F_n)$, the strong Markov property applies at the $(\F_n)$-stopping times $\rho_r$, and the cycle segments are therefore i.i.d.

It remains to show that the regeneration times are almost surely finite. Let
\[
   C_K:=\{y:\ V(y)\le K\}
\]
with $K$ as in Theorem~\ref{theo:poisson-stability}, and let
\[
   \sigma_0:=\inf\{n\ge 0:\ \mathsf{Y}_n\in C_K\},\qquad
   \sigma_{j+1}:=\inf\{n>\sigma_j:\ \mathsf{Y}_n\in C_K\},\quad j\ge 0,
\]
be the successive return times to $C_K$. By Theorem~\ref{theo:poisson-stability}, the chain $(\mathsf{Y}_n)$ is positive Harris recurrent, so all $\sigma_j$ are almost surely finite. By Proposition~\ref{prop:small-set-poisson} with $M=K$, whenever $\mathsf{Y}_{\sigma_j}\in C_K$, the event of $m_K+2$ consecutive long gaps after time $\sigma_j$ has conditional probability $p_K>0$ and forces a regeneration at time $\sigma_j+m_K+2$. Hence
\[
   \P\bigl(\rho_0>\sigma_j+m_K+2 \,\big\vert\, \F_{\sigma_j}\bigr)\le 1-p_K
   \qquad\text{a.s.\ on }\{\rho_0>\sigma_j\}.
\]
Taking expectations and using induction over $j$ yields
\[
   \P(\rho_0>\sigma_j+m_K+2)\le (1-p_K)^{j+1}\longrightarrow 0,
\]
so $\rho_0<\infty$ almost surely. The same argument applied after the stopping time $\rho_r$ shows that $\rho_{r+1}<\infty$ almost surely whenever $\rho_r<\infty$. Hence all regeneration times are almost surely finite.
\end{proof}

\begin{cor}[Stationary limits and ergodic averages]\label{cor:regenerative-ergodic}
There exist integrable random variables $X_\infty^{(o)}$ and $X_\infty^{(i)}$, $i\in\Nodes$, such that
\[
   X_n^{(o)} \Longrightarrow X_\infty^{(o)},
   \qquad
   X_n^{(i)} \Longrightarrow X_\infty^{(i)},
   \qquad i\in\Nodes,
\]
and, with $\mu^{(o)}:=\E[X_\infty^{(o)}]$ and $\mu^{(i)}:=\E[X_\infty^{(i)}]$,
\[
   \frac1n\sum_{m=1}^n X_m^{(o)} \longrightarrow \mu^{(o)},
   \qquad
   \frac1n\sum_{m=1}^n X_m^{(i)} \longrightarrow \mu^{(i)},
   \qquad i\in\Nodes,
\]
almost surely.
\end{cor}

\begin{proof}
Let $\pi$ denote the invariant probability measure from Theorem~\ref{theo:poisson-stability}. Since $(\mathsf{Y}_n)$ is positive Harris recurrent and aperiodic, Meyn and Tweedie~\cite[Theorem~13.3.3]{MeynTweedie2009} imply that
\[
   \|\mathcal{L}(\mathsf{Y}_n)-\pi\|_{\mathrm{TV}}\to 0
\]
for every initial distribution. Therefore $X_n^{(o)}$ and $X_n^{(i)}$ converge in distribution to the corresponding $\pi$-marginals, which we denote by $X_\infty^{(o)}$ and $X_\infty^{(i)}$.

To prove integrability, start the chain from the deterministic genesis state, for which $V_0=1\le 1+\alpha$. The exponential bound \eqref{eq:exp-moment} with $m=0$ then gives
\[
   \E[e^{\eta V_n}]
   \le e^{\eta} + \frac{1-\rho^n}{1-\rho}\,e^{\eta(1+\alpha)}\,e^{\eta K}
   \le C_\eta,
   \qquad n\ge 0,
\]
for the finite constant
\[
   C_\eta:= e^{\eta} + \frac{e^{\eta(1+\alpha)}e^{\eta K}}{1-\rho}.
\]
For $M>0$ define the bounded measurable truncation
\[
   f_M(y):=e^{\eta V(y)}\wedge M.
\]
By total-variation convergence,
\[
   \pi(f_M)=\lim_{n\to\infty}\E[f_M(\mathsf{Y}_n)]\le C_\eta.
\]
Letting $M\to\infty$ and using monotone convergence yields
\[
   \pi(e^{\eta V})\le C_\eta<\infty.
\]
In particular, $\pi(V)<\infty$. Since
\[
   X_n^{(o)}\le V_n,
   \qquad
   X_n^{(i)}\le X_n^{(o)}+kL_n\le \Bigl(1+\frac{k}{\alpha}\Bigr)V_n,
\]
the stationary random variables $X_\infty^{(o)}$ and $X_\infty^{(i)}$ are integrable.

Finally, Meyn and Tweedie~\cite[Theorem~17.1.7]{MeynTweedie2009} state that for a positive Harris chain the strong law of large numbers holds for every function in $L^1(\pi)$. Applying this to the measurable coordinate functions
\[
   f_o(y):=X^{(o)}(y),
   \qquad
   f_i(y):=X^{(i)}(y),
\]
gives the almost-sure convergence of the empirical averages to their stationary means.
\end{proof}

\begin{rem}
The regeneration structure above is qualitative: it identifies i.i.d.\ cycle segments for the embedded chain. The quantitative exponential bounds in Section~\ref{sec:quantitative} are obtained instead from a direct verification of Hajek's exponential-drift hypotheses for the Lyapunov function $V_n$.
\end{rem}

\begin{rem}[Relation to the hard-core regeneration argument]
The hard-core case is technically simpler because bounded inter-arrival gaps imply bounded recent activity. Then the discrepancy between observer, common, and local tip pools is uniformly controlled, so once $X_n^{(o)}$ is large the next issuer sees mostly common tips and a negative drift estimate follows directly. In the Poisson setting the number of recent arrivals is unbounded, so one must control the extra fluctuation carried by $L_n$, which is why the augmented Lyapunov function $V_n=X_n^{(o)}+\alpha L_n$ is needed.
\end{rem}

\section{A Little-type identity}\label{sec:little}
We now formulate Little's law only in the stationary Poisson regime. Let $\pi$ denote the invariant law of the embedded chain from Theorem~\ref{theo:poisson-stability}. For each issued block $m\ge1$, define its \emph{observer tip lifetime} by
\[
   W_m := \inf\{t_n:\ n>m,\ \blockID_m \in \refBlock_n\}-t_m \;\in\; [0,\infty],
\]
with the convention $\inf\emptyset=\infty$.

\begin{prop}[Stationary marked-point-process representation]\label{prop:little-palm}
There exists a stationary marked point process
\[
   \Psi=\sum_{n\in\mathbb Z}\delta_{(T_n,W_n)}
\]
on $\mathbb{R}\times [0,\infty]$ with intensity $\lambda$, where $W_n$ is the observer-tip lifetime of the block issued at time $T_n$, whose Palm version $\P^0$ satisfies $T_0=0$ almost surely and has post-issuance embedded state at time $0$ distributed according to $\pi$. Under $\P^0$, the next inter-arrival time $T_1-T_0$ is $\mathrm{Exp}(\lambda)$ and is independent of the post-issuance state at time $0$.
\end{prop}

\begin{proof}
Because Theorem~\ref{theo:poisson-stability} yields an invariant law $\pi$ for the enlarged embedded chain, there exists a two-sided stationary version $(\mathsf{Y}_n)_{n\in\mathbb Z}$ with one-dimensional marginal $\pi$. Realise each transition by an innovation consisting of the next exponential inter-arrival gap, issuer identity, block identifier, delay vector, and reference-selection randomness. The gap component is $\mathrm{Exp}(\lambda)$ and is independent of the current post-issuance state, as in Assumption~\ref{ass:poisson}. These innovations determine a two-sided stationary issuance-by-issuance evolution and therefore the Palm version of the marked issuance process with $T_0=0$. The lifetime mark $W_n$ is then defined from this bi-infinite evolution exactly as in the one-sided model. The Palm inversion theorem for stationary point processes yields the associated time-stationary marked point process $\Psi$; see \cite[Chapter~VI]{Asm:03}. Its intensity is $\lambda$, because under the Palm version the successive inter-arrival times are $\mathrm{Exp}(\lambda)$.
\end{proof}

\begin{cor}[Little-type identity under Poisson stability]\label{cor:little-poisson}
Let $\P^0$ denote the Palm probability from Proposition~\ref{prop:little-palm}, and let $W_0$ be the lifetime mark of the block issued at time $0$. Then
\[
   \pi(X^{(o)}) := \int X^{(o)}(y)\,\pi(dy) = \lambda\,\E^{0}[W_0].
\]
Equivalently, if $(\bar X_t^{(o)})_{t\in\mathbb R}$ denotes the time-stationary observer tip-count process associated with $\Psi$, then
\[
   \E[\bar X_0^{(o)}] = \lambda\,\E^{0}[W_0] = \pi(X^{(o)}).
\]
In particular, $\E^{0}[W_0]<\infty$.
\end{cor}

\begin{proof}
Let $(\bar X_t^{(o)})_{t\in\mathbb R}$ be the observer tip-count process generated by the stationary marked point process $\Psi$. A block issued at time $u$ contributes to the observer tip count during the half-open lifetime interval $[u,u+w)$. Pathwise, the number of observer tips present at time $0$ is therefore the number of lifetime intervals covering $0$, namely
\[
   \bar X_0^{(o)}
   = \int \1\{u\le 0 < u+w\}\,\Psi(du,dw).
\]
Applying Campbell's formula under the stationary law of $\Psi$ gives
\[
   \E[\bar X_0^{(o)}]
   = \lambda\,\E^{0}\!\left[\int_{\mathbb R}\1\{u\le 0 < u+W_0\}\,du\right]
   = \lambda\,\E^{0}[W_0].
\]

To relate this time-stationary mean to the embedded stationary mean, apply Palm inversion once more:
\[
   \E[\bar X_0^{(o)}]
   = \lambda\,\E^{0}\!\left[\int_0^{T_1} \bar X_t^{(o)}\,dt\right].
\]
Under $\P^0$, the observer tip count changes only at issuance times, so $\bar X_t^{(o)}=X_0^{(o)}$ for $0\le t<T_1$, where $X_0^{(o)}$ denotes the post-issuance observer tip count at time $0$. Hence
\[
   \E[\bar X_0^{(o)}]
   = \lambda\,\E^{0}[T_1 X_0^{(o)}]
   = \lambda\,\E^{0}[T_1]\,\E^{0}[X_0^{(o)}]
   = \pi(X^{(o)}),
\]
because $T_1\sim\Exp(\lambda)$ is independent of $X_0^{(o)}$, $\lambda\,\E^{0}[T_1]=1$, and $X_0^{(o)}$ has law $\pi$ under $\P^0$ by Proposition~\ref{prop:little-palm}. Combining the two identities yields
\[
   \pi(X^{(o)}) = \lambda\,\E^{0}[W_0].
\]
Finally, Corollary~\ref{cor:regenerative-ergodic} gives $\pi(X^{(o)})<\infty$, hence $\E^{0}[W_0]<\infty$ as well.
\end{proof}

\begin{rem}
Corollary~\ref{cor:little-poisson} is the standard ``mean number in system equals arrival rate times mean sojourn time'' identity, but the proof must distinguish carefully between the embedded stationary law $\pi$ at issuance epochs and the time-stationary Poisson regime obtained from it by Palm inversion. In the present model these means coincide for the observer tip count because the observer process is constant between issuances and the Poisson holding times are state-independent.
\end{rem}

\section{Illustrative simulations}\label{sec:simulations}
We provide illustrative simulation results to complement the theory. Their role is not to establish quantitative performance estimates or to validate the stability proof, but to show qualitatively how local tip-pool perceptions vary with delay randomness and issuance heterogeneity in more detailed network settings where explicit formulas are unavailable. The plotted trajectories include the initial transient from the genesis state and should therefore be read as protocol illustrations rather than as steady-state estimators.
The simulations are performed in an open-source simulator \cite{otv-simulator} also used in \cite{robustnessIOTA2.0}. This simulator models both communication over the peer-to-peer layer and block creation.
The statistical analysis of the data is done with the software R (4.1.2), and the package ``ggstatsplot'' \cite{ggstatsplot}. 

We use a gossip protocol to model the network latency on a network topology with a small diameter. More precisely, we use a Watts-Strogatz network \cite{Watts1998collective} with mean degree $10$ and re-wiring probability $1$. 
The gossip algorithm forwards new blocks to all neighbours in the Watts-Strogatz network. The delay for each connection on the peer-to-peer layer is independent and uniformly distributed in the interval $[\delta_{min}, \delta_{max}]$.

We model the different issuance rates of the nodes in the network using the Zipf empirical law with parameter $s$~\cite{zipfs-law}. This is motivated by the fact that in a real-world scenario with heterogeneous weights, the Zipf law is frequently observed, e.g., see~\cite{Adamic2002ZipfsLA,wealth_pareto,Li2002ZipfsLE}. Note that with Zipf's law, a homogeneous network, e.g., can be modelled for $s=0$, while the higher the $s$, the more heterogeneous or centralized the weight distribution becomes. 

\subsection{Heterogeneous rates}
The issuing rates of the $N=100$ nodes are Zipf-distributed with parameter $s$, i.e.,
\begin{equation}
    \lambda_i = \frac{i^{-s}}{\sum_{j=1}^N j^{-s}} \lambda,
\end{equation}
where $\lambda$ is the total issuance rate. 

We have set the other parameters of our numerical experiments as follows: the number of references $k=8$.  This choice of $k=8$ is made since it is in the ``middle'' on a logarithmic scale of the extreme cases $2^0$ and $2^7$.  If $k=1$ we obtain a tree and if $k$ is close to the number of nodes, then the number of tips is generally very small. Moreover, $k=8$ is the value considered in \cite{robustnessIOTA2.0}.

The network latency between two peers in the peer-to-peer network is modelled by a uniform random variable with
$\delta_{min}=20ms,$ $\delta_{max}=180ms$. It is a common assumption to consider the mean latency to be close to $100$ms. Moreover, most delays in wide area networks and the Internet fall into our interval, e.g., see \cite{internetLatency}.
The total block issuance rate is set to $\lambda=500$ blocks per second (BPS). The local tip pools are measured in the simulation every $50ms$, and every simulation lasts for $60$ seconds.

Let us first consider the case of a heterogeneous node activity, $s=1$. In this scenario, Node~$1$ issues blocks at a rate of $96$ BPS, Node~$2$ with a rate of $48$ BPS, and the two slower nodes, Node $99$ and $100$ issue with rates around $1$ BPS.

Figures~\ref{subfig:sim1} and~\ref{subfig:sim1b} show the resulting local tip-pool trajectories and empirical sample distributions for these nodes.
\subsection{Homogeneous rates}
We consider the homogeneous case, where every node issues blocks with the same rate, i.e. $s=0$.  The other parameters are set as before. The results in Figures~\ref{subfig:sim2} and~\ref{subfig:sim2b} show that the local tip pools have similar sizes in this illustrative run. Comparing these results with the heterogeneous setting in Figure~\ref{subfig:sim1b} suggests that stronger concentration of issuance activity can reduce the observed tip-pool sizes, although the experiment is not designed to provide a quantitative estimate of this effect.

\subsection{Effect of delay randomness}
\begin{figure}
     \centering
     \begin{subfigure}[b]{\textwidth}
         \centering
         \includegraphics[width=0.78\textwidth]{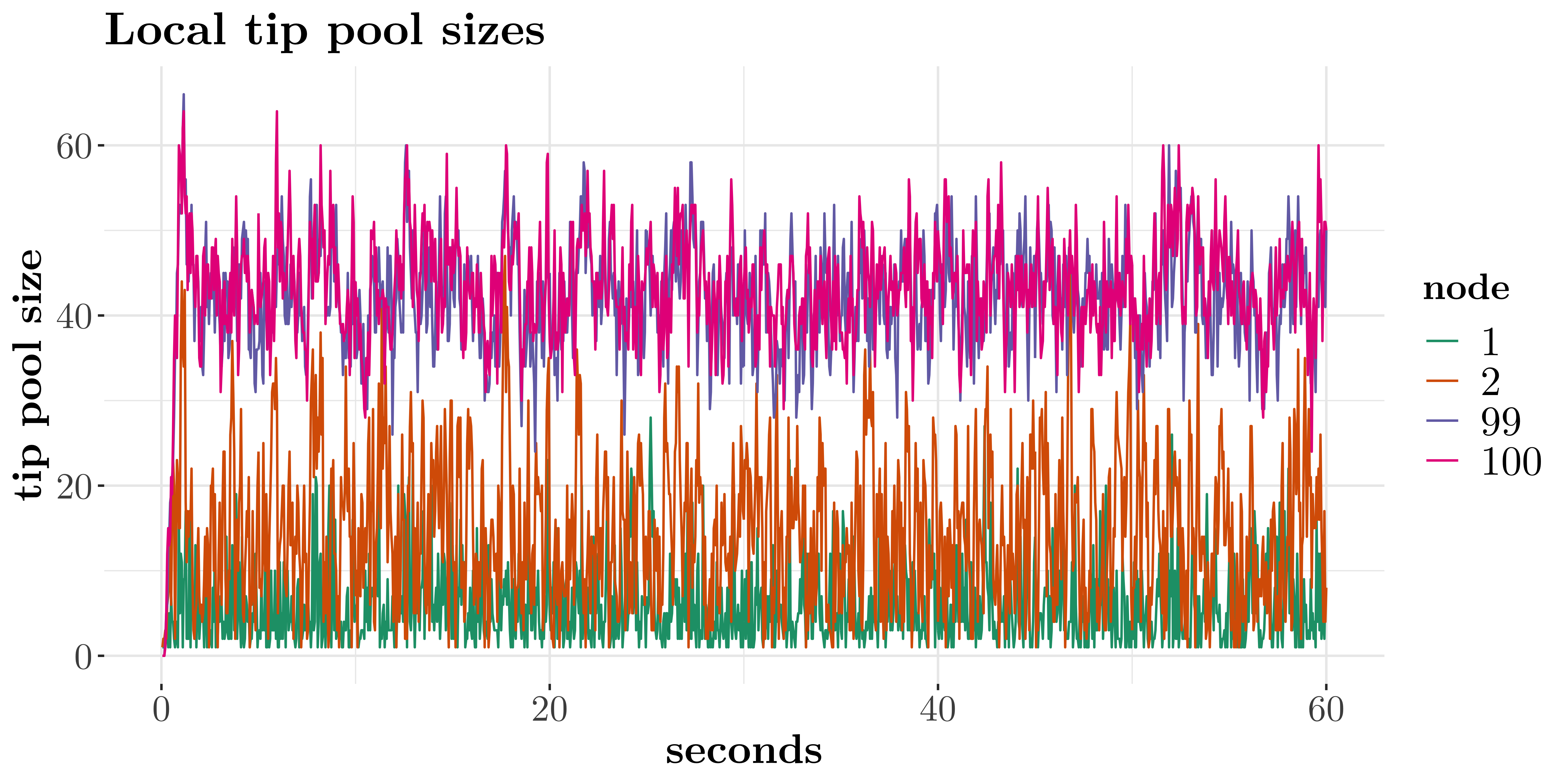}
         \caption{Heterogeneous rates according to Zipf law with $s=1$, BPS of $500$, and random network delay.}
         \label{subfig:sim1}
     \end{subfigure}\\
     \begin{subfigure}[b]{\textwidth}
         \centering
         \includegraphics[width=0.78\textwidth]{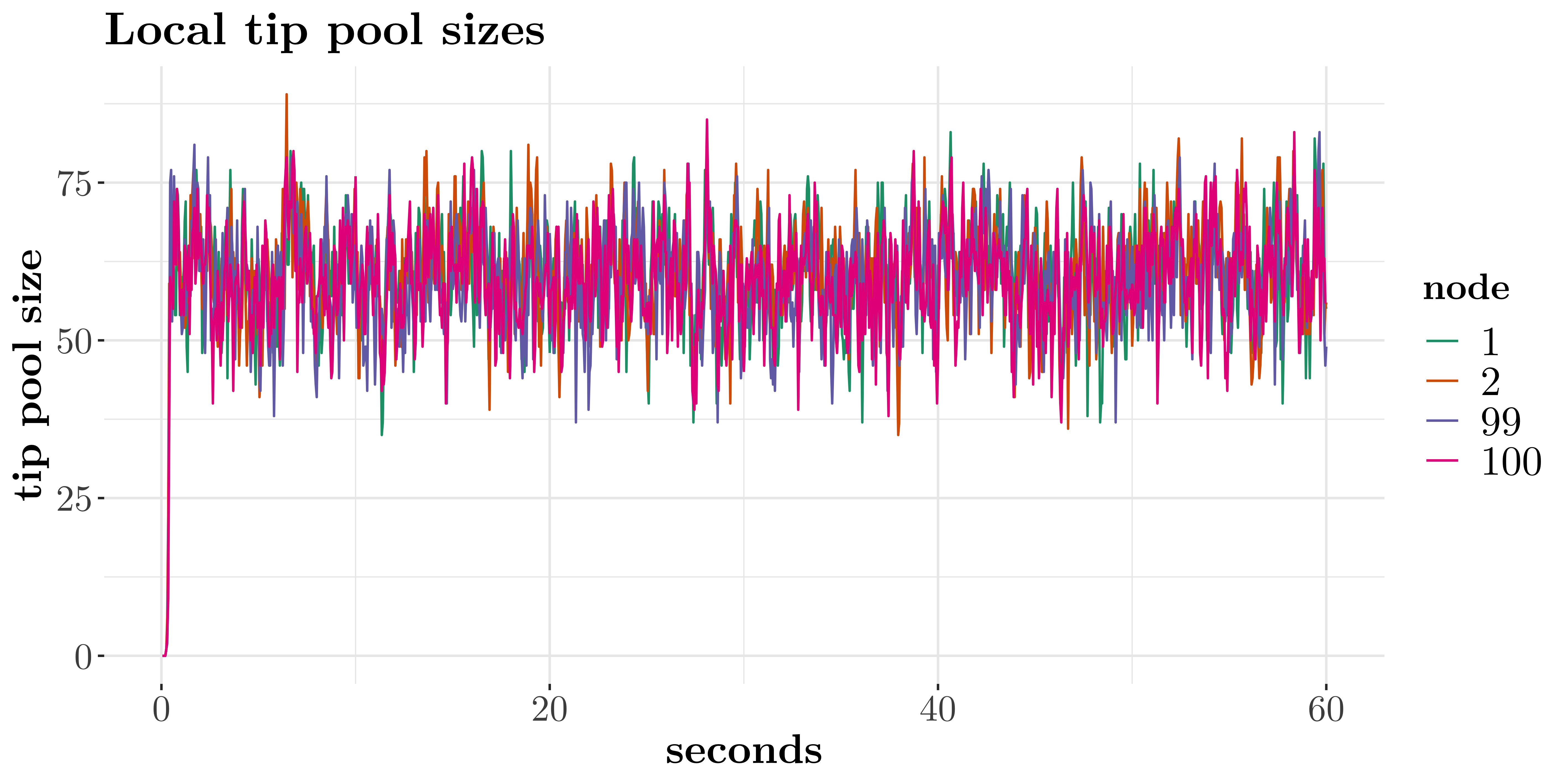}
         \caption{Homogeneous rates according to Zipf law with $s=0$, BPS of $500$, and random network delay.}
         \label{subfig:sim2}
     \end{subfigure}\\
     \begin{subfigure}[b]{\textwidth}
         \centering
         \includegraphics[width=0.78\textwidth]{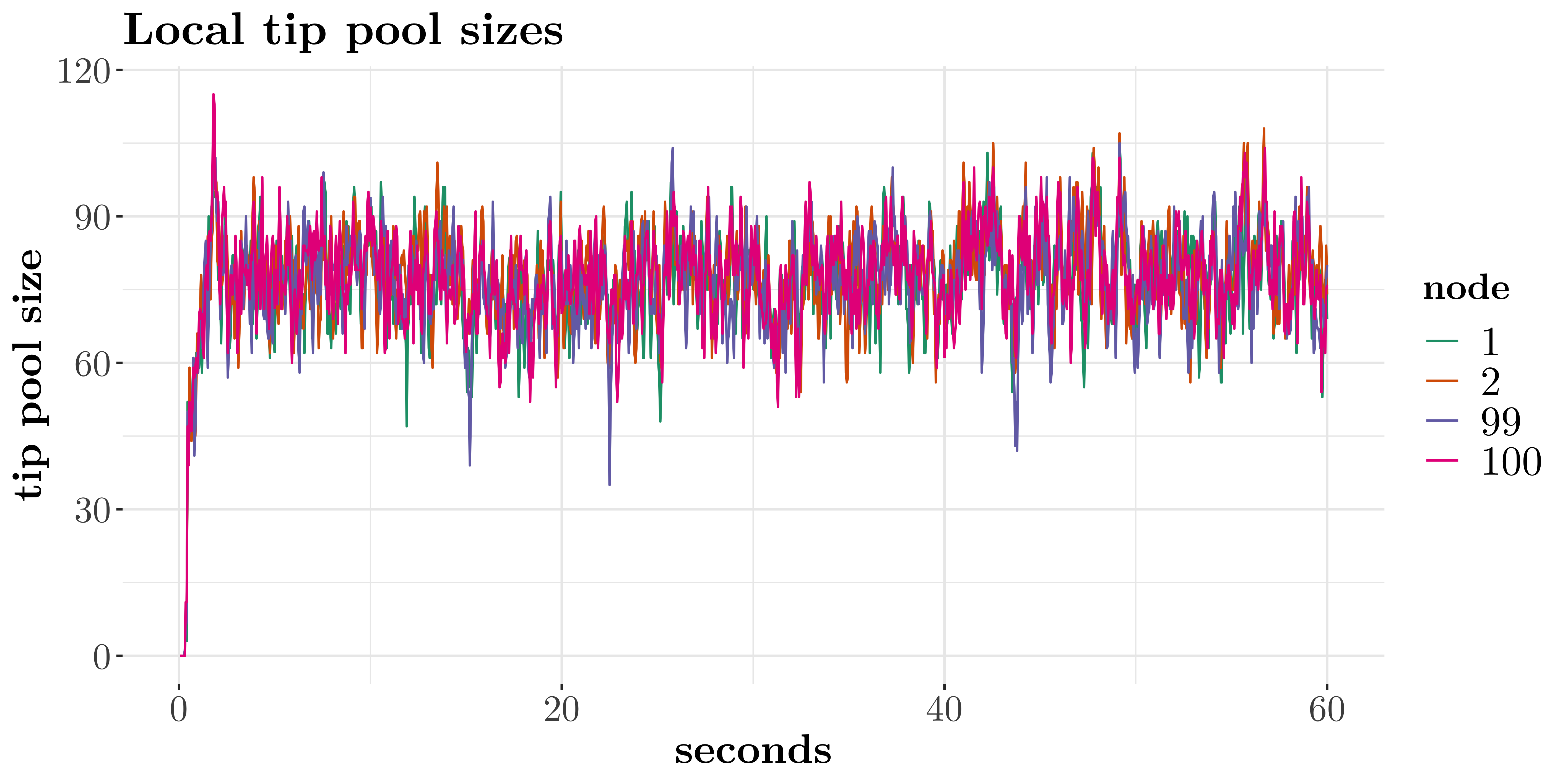}
         \caption{Homogeneous rates according to Zipf law with $s=0$, BPS of $500$, and constant network delay of $100$ms.}
         \label{subfig:sim3}
     \end{subfigure}
        \caption{Illustrative tip-pool trajectories of selected high- and low-rate nodes with $N=100$ nodes under three scenarios. The trajectories include the initial transient from the genesis state.}
        \label{fig:sims}
\end{figure}

\begin{figure}
     \centering
     \begin{subfigure}[b]{\textwidth}
         \centering
         \includegraphics[width=0.76\textwidth]{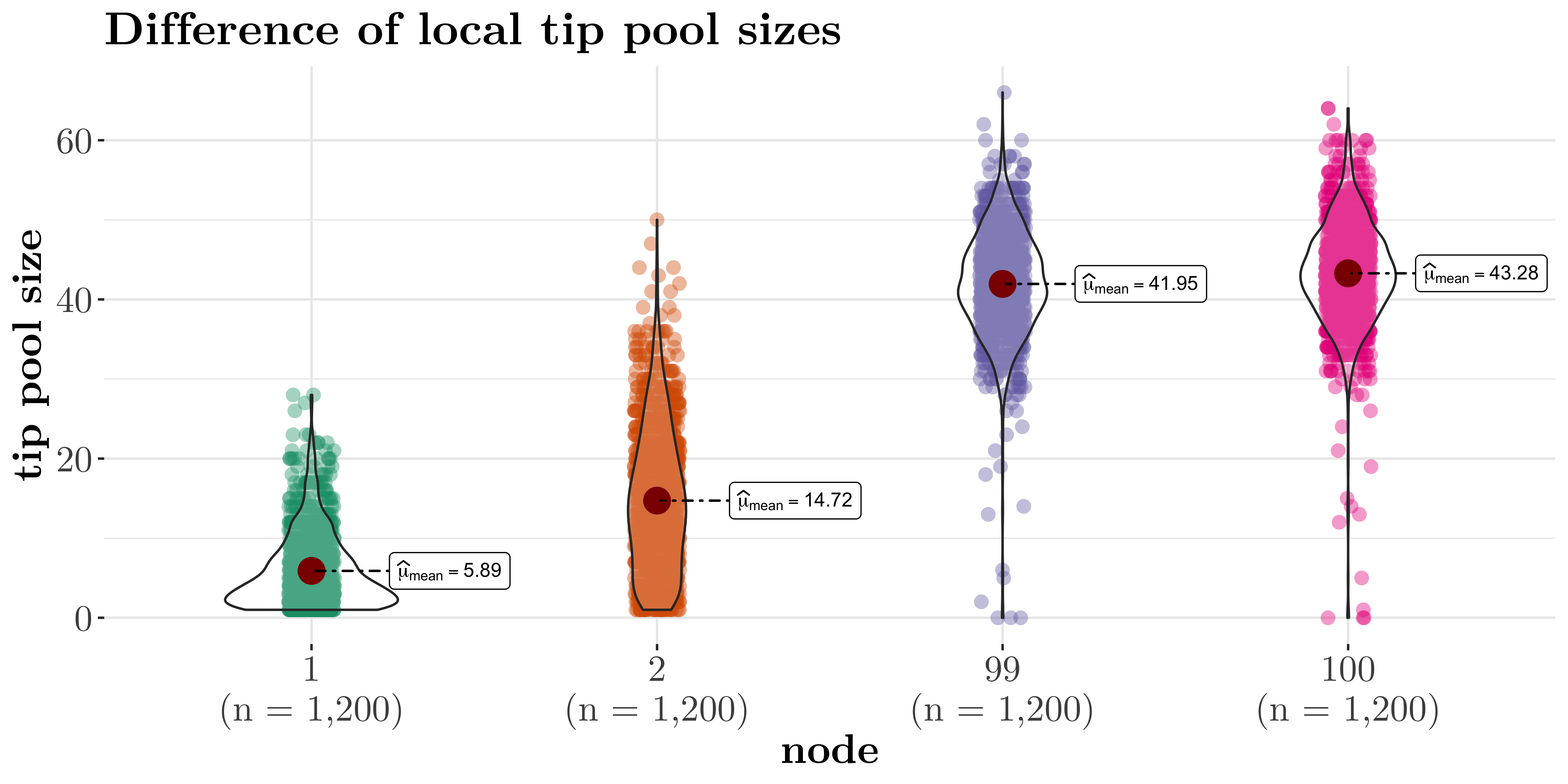}
         \caption{Heterogeneous rates according to Zipf law with $s=1$, BPS of $500$, and random network delay.}
         \label{subfig:sim1b}
     \end{subfigure}\\
     \begin{subfigure}[b]{\textwidth}
         \centering
         \includegraphics[width=0.76\textwidth]{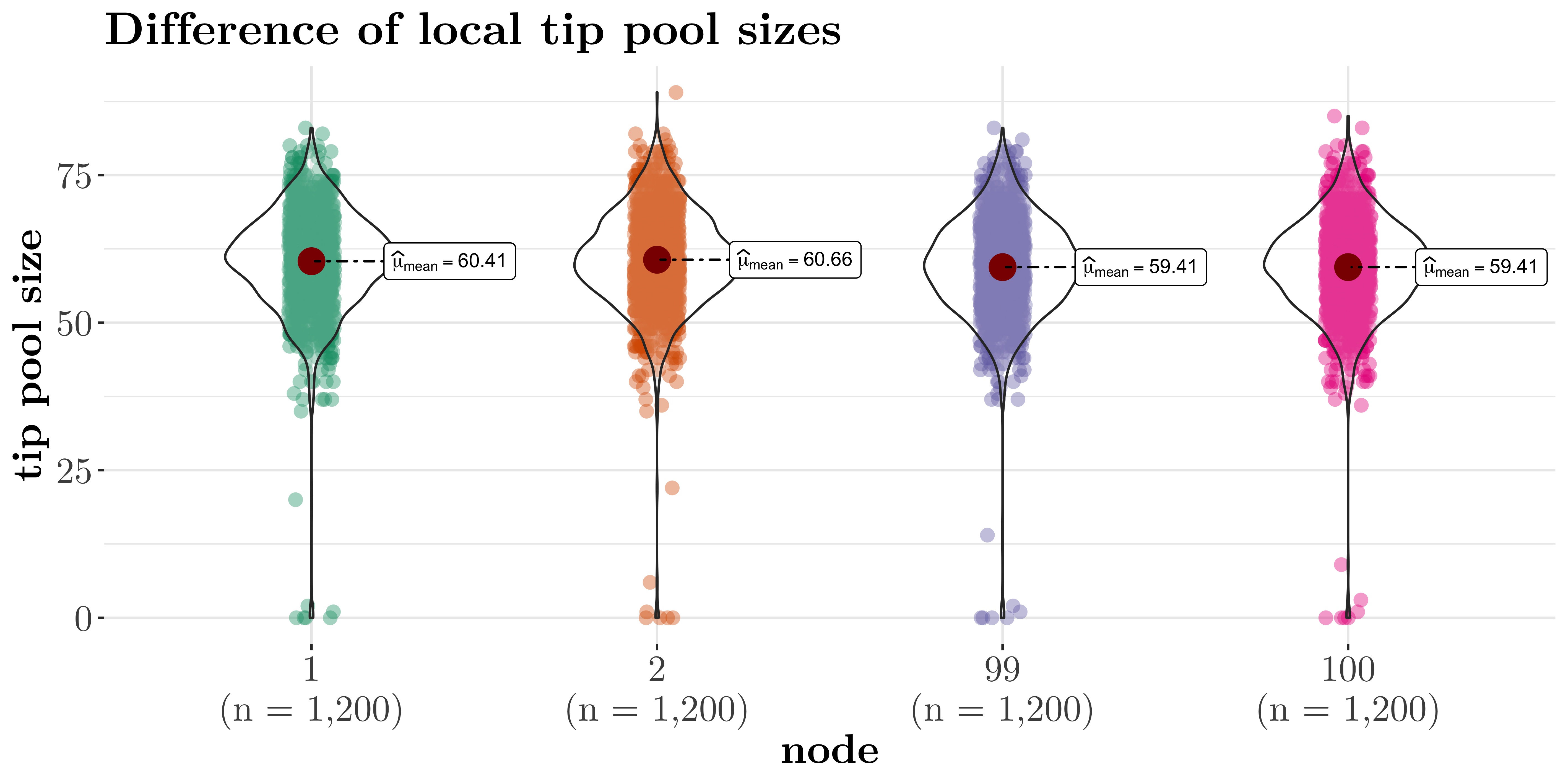}
         \caption{Homogeneous rates according to Zipf law with $s=0$, BPS of $500$, and random network delay.}
         \label{subfig:sim2b}
     \end{subfigure}\\
     \begin{subfigure}[b]{\textwidth}
         \centering
         \includegraphics[width=0.76\textwidth]{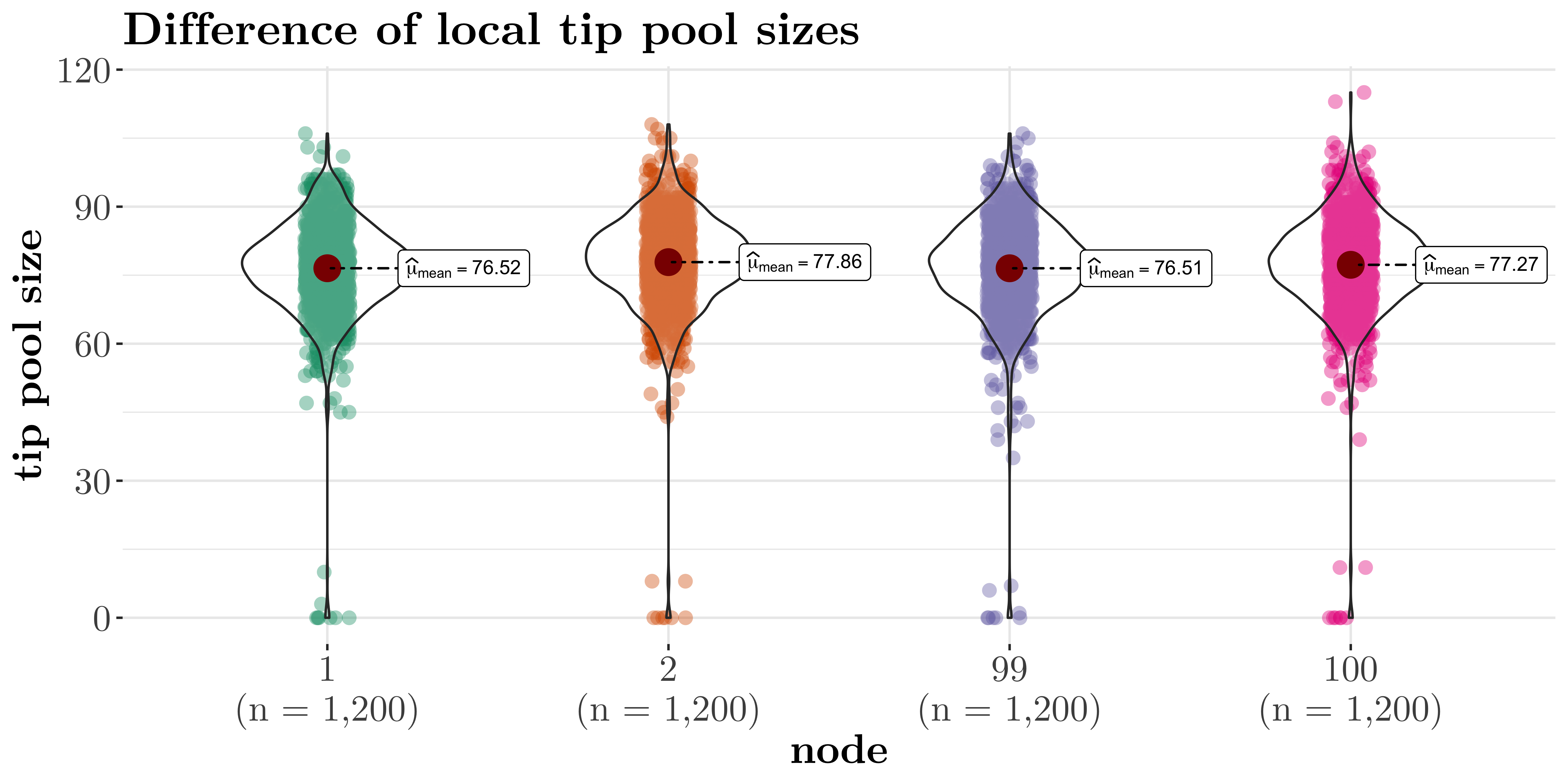}
         \caption{Homogeneous rates according to Zipf law with $s=0$, BPS of $500$, and constant network delay of $100$ms.}
         \label{subfig:sim3b}
     \end{subfigure}
        \caption{Illustrative empirical distributions of sampled local tip-pool sizes with $N=100$ nodes under the three scenarios of Figure~\protect\ref{fig:sims}. The samples are pooled from the plotted simulation runs and are not intended as steady-state confidence intervals.}
        \label{fig:simsb}
\end{figure}

The preceding illustrative scenarios suggest that different issuing rates can considerably affect the local tip pools. A natural explanation is that the average delay of high-frequency nodes is much smaller than those of lower frequencies. Previous heuristic results~\cite{kumar2022effect} also indicate that the full delay distribution, not only its mean, can affect tip-pool sizes. We illustrate this effect by comparing the homogeneous-rate case $s=0$ with uniformly distributed edge delays on $[20,180]$ms against the same setup with constant $100$ms delays, see Figures~\ref{subfig:sim3} and~\ref{subfig:sim3b}. In this illustrative run, the constant-delay scenario shows larger tip pools than the uniform-delay scenario with the same mean delay. This effect is also present for heterogeneous rates, but we omit the figures for brevity.

\section{Discussion and extensions}\label{sec:discussion}
This paper presents a continuous-time model for DAG-based distributed ledgers with variable and heterogeneous network delays. We assume partial synchrony (Assumption~\ref{ass:sync}) and global Poisson issuance (Assumption~\ref{ass:poisson}). Under these assumptions we established stability and a stationary regime for the observer and local tip pool sizes via the augmented Lyapunov argument of Section~\ref{sec:poisson}, obtained exponential tail and moment bounds through the Hajek analysis of Section~\ref{sec:quantitative}, derived a Little-type identity in the stationary Poisson regime in Section~\ref{sec:little}, and complemented the theory by illustrative Monte-Carlo simulations. We close with several directions for further work.

\textbf{Different types of delays:}
As already mentioned in Subsection~\ref{sec:relworks}, a different type of delay, namely the time to validate a block, has been studied in~\cite{penzkofer2021impact}.  One natural way to incorporate such delays is to include an additional mark in the Poisson point process that encodes the block type. The delays of a block then also depend on its type. 
The proof strategy should extend to this more general situation, but understanding how such delays affect tip-pool sizes would require more quantitative estimates.

\textbf{Quantitative results} 
Beyond qualitative stability, we obtained exponential tail and return-time bounds for the Lyapunov function and therefore polynomial moment bounds for the observer and local tip pool sizes. For applications, however, one would like sharper explicit performance bounds directly in terms of the delay distribution and the issuance heterogeneity. The most important quantity is the expected tip pool size. Previous results, \cite{kumar2022effect,penzkofer2021impact, cullen2019-variable-delay-9142287}, and our simulations show that the tip pool size depends strongly on the full delay distribution, so closed-form formulas seem out of reach in general. A more feasible goal is to derive robust upper and lower bounds that capture this dependence more explicitly and quantify the speed of convergence to stationarity.

\textbf{Extreme values and large deviations: } A natural next step is to sharpen the exponential tail bounds of Section~\ref{sec:quantitative} into genuine large-deviation and extreme-value asymptotics for the tip pool process. Such estimates would quantify the risk of atypically large tip pools more precisely and would be directly relevant for protocol dimensioning. Our drift and regeneration analysis suggests that this should be approachable once sharper quantitative return-time estimates are available.

\subsection{Hard-core issuance as a simpler variant}\label{sec:hardcore-variant}
If one imposes a hard-core lower bound $\delta_\ast>0$ on the global inter-arrival times, then every interval of length $\Delta$ contains at most $B=\lfloor \Delta/\delta_\ast\rfloor+1$ issuance events. Hence the recent-activity count is uniformly bounded, and Lemma~\ref{lem:X_nX_nc} yields a uniform bound between observer, common, and local tip pools. For large observer tip pools, almost all visible tips are then common, so the next issuer references common tips with high probability and the negative-drift argument becomes simpler than in the Poisson case, where the unbounded window count must be controlled explicitly through the Lyapunov term $\alpha L_n$.

\bibliographystyle{plain}
\bibliography{references}

\end{document}